%% file: FineRapport.tex
\documentclass[10pt,a4paper]{article}

\usepackage[latin1]{inputenc} 
\usepackage{graphicx}
\usepackage{epsf}
\usepackage{amsmath}
\usepackage{multirow}
\usepackage[usenames]{color}
\usepackage{gastex}

\newtheorem{theorem}{Theorem}[section] 
\newtheorem{lemma}[theorem]{Lemma}
\newtheorem{property}[theorem]{Property}

\newtheorem{definition}[theorem]{Definition}

\newenvironment{proof}[1][Proof]{\begin{trivlist} 
\item[\hskip \labelsep {\bfseries #1}]}{$ $ \qed \end{trivlist}}
\newenvironment{example}[1][Example ]{\begin{trivlist} 
\item[\hskip \labelsep {\bfseries #1}]}{\end{trivlist}} 
\newenvironment{remark}[1][Remark ]{\begin{trivlist} 
\item[\hskip \labelsep {\bfseries #1}]}{\end{trivlist}}

\newcommand{\qed}{\sc q.e.d.}

\input{Commd_bibl}

\newcommand{\ffineRogers}[1]{\frac{1}{2} \sum_{i=0}^{#1-1} (-\frac{1}{2})^{i} C_{#1+1-i}}
\newcommand{\ffineStrehl}[1]{\sum_{k=1}^{\lfloor\frac{#1+1}{2}\rfloor} C_{n+1,2k}}
\newcommand{\fcatalan}[1]{\frac{1}{#1 + 1} {{2 #1} \choose #1}}

\newcommand{\fcatalanDelannoy}[2]{{{2#1-#2-1} \choose {#1-1}} - {{2#1-#2-1} \choose {#1}}}

\def\bbbn{{\rm I\!N }}

\author{O.~Guibert and S.~Pelat-Alloin}

\title{Extending Fine sequences: a link with forbidden patterns}

\begin{document}

\maketitle

\begin{abstract}
We propose a natural, bivariate, generalization of the nonsingular 
similarity relations considered by T.~Fine. We also provide an
enumeration formulae and a generating tree for those relations. The
latter allow us to give a new bijection between
$321$-avoiding derangements and Fine sequences. Moreover, we establish that
two special cases are in a one-to-one correspondence with subsets of
permutations characterized by forbidden subsequences on the
symmetrical group. All our results are established using the technique
of generating tree, thus giving entirely bijective proofs.
\end{abstract}


\bigskip

\section{Introduction and motivation}

\subsection{Similarity relations and permutations with forbidden patterns}
\label{sr}

\emph{Similarity relations} are \emph{symmetrical} and
\emph{reflexive} binary relations $R$ operating on
$[n]=\{1,2,\ldots,n\}$ such that $x < y < z$ and $xRz$ implies $xRy$
and $yRz$. These relations can be coded by an \emph{integer sequence}
$\alpha_1 \alpha_2 \ldots \alpha_n$ such that for every $y$:
$\alpha_y = y - x$ where $x$ is the minimal integer verifying
$xRy$. From the above definitions, one can see that $\alpha_1=0$ and
$0 \leq \alpha_{x+1} \leq \alpha_x + 1$. Such relations are enumerated
by the Catalan numbers $C_n = \fcatalan n$ for $n \geq 0$. The figure
below illustrates the similarity relation on $n=8$ coded by the
integer sequence $01100121$ (an edge connects two vertices $x,y$ iif
$xRy$).

\begin{figure}[h]
\center
\vspace*{-5.5truecm}
{\centerline{\epsfxsize=300pt\epsfbox{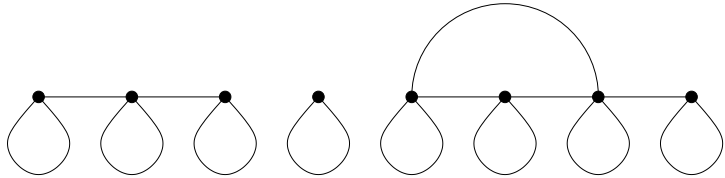}}}
\vspace*{-5.5truecm}
\caption{a similarity relation on $n=8$ coded by the
integer sequence $01100121$.}
\end{figure}

In the present work, we shall consider \emph{nonsingular similarity
  relations} on $[n]$, that is similarity relations with the
restriction that every element must be in relation with another one at
least, namely: $\forall x \in [n], \exists y \not = x$ such that
$xRy$. We can code these relations on $[n]$ with the same integer
sequence as above by adding the condition that every $0$ is followed
by $1$. From now on, the cardinality of these relations on $[n]$ will
be denoted $F_n$, with $n \geq 1$, the first values of which are:
$1,2,6,18,57,186,622,2120,\ldots$ These are the Fine sequences,
enumerated by $F_n = \ffineRogers n$, a formulae due to D.G.~Roger
\cite{Rogers}. The figure hereafter illustrates a nonsingular
similarity relation on $n=6$ coded by the integer sequence $010122$.

\begin{figure}[h]
\center
\vspace*{-5.5truecm}
{\centerline{\epsfxsize=300pt\epsfbox{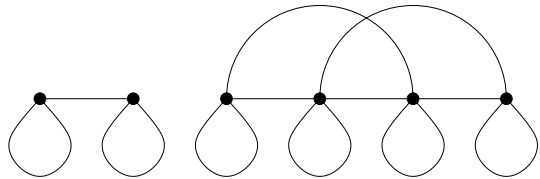}}}
\vspace*{-5.5truecm}
\caption{a nonsingular similarity relation on $n=6$ coded by the integer sequence $010122$.}
\end{figure}

This Fine sequence was considered notably by L.W.~Shapiro
\cite{Shapiro} who stated that $2 F_n + F_{n-1} = C_{n+1}$. 
It was also discussed by V.~Strehl \cite{Strehl} who established that
$F_n = \ffineStrehl{n}$ for $n \geq 1$ 
where $C_{n,j} = \fcatalanDelannoy{n}{j}$ enumerates the nonsingular
similarity relation on $[n+1]$ for which the transitive 
closure consists of $k$ blocks. We recall that $C_{n,j}$ is known as ballot numbers, Delannoy numbers \cite{Errera} or 
distribution $\alpha$ of the Catalan numbers \cite{KrewEventail}. More
recently, Deutsch and Shapiro wrote a survey \cite{Deutsch} and Callan give also some
identities for the Fine numbers \cite{Callan}.

The study of \emph{permutations with forbidden patterns} can be
 traced back to Simion and Schmidt \cite{Simion}. It is, nowadays, a
 growing domain of combinatorics with its own annual conference (for
 instance see \cite{Elder}). 

Permutations with forbidden patterns (or subsequences) constitute a subset of the
symmetrical group characterized by the exclusion of permutations
containing at least one subsequence orderisomorphic to the forbidden
one. More specifically, a permutation $\pi$ of length $n$ contains the
subsequence (pattern) of type $\tau$ of length $k$ if and only if
one can find $1 \leq i_{\tau(1)} < i_{\tau(2)} < \cdots < i_{\tau(k)}
\leq n$ such that $\pi(i_1) < \pi(i_2) < \cdots < \pi(i_k)$. We denote
by ${\cal S}_n (\tau)$ the set of permutations of length $n$ which
does not contain any subsequences of type $\tau$. Moreover, let ${\cal
  S}_n (\tau_1,\tau_2,\ldots,\tau_l )$ be ${\cal S}_n(\tau_1) \cap
{\cal S}_n(\tau_2) \cap \ldots \cap {\cal S}_n(\tau_l)$. As an example, the
following permutation $\pi = 364125$ contains two occurences of the
pattern $123$ (namely the subsequences $345$ and $125$) but none of the pattern
$321$, so it belongs, for example, to ${\cal S}_6 (321)$. \\

The aim of our present work is to give a natural generalization of
nonsingular similarity relations. In section \ref{defgfs}, we define
these generalized Fine sequences as both Catalan paths with
constraints or words on the natural numbers. We also give an
enumeration formulae for those objects. Thereafter, in section
\ref{mt}, we state four theorems. The first one relates the
generalized Fine sequences with generating trees. The second one
establish a link between those generating trees and $321$-avoiding derangements, which ones are notably in
bijection with non singular similarity relation \cite{Deutsch}. The others relate
two cases of generalized Fine sequences with permutations with
forbidden patterns. Next, in section \ref{bk}, we introduce some
needed background on those topics. We then proceed to demonstrate
these theorems, the former in section \ref{ptgfs} and the latter ones
in sections \ref{sder}, \ref{pfsf} and \ref{pfsh}.

\section{Generalized Fine sequences}
\label{defgfs}

A natural way of generalizing nonsingular similarity relations is to
consider a Fine sequence as a Catalan path and operate a congruous/modulo fonction 
on the $\alpha$-distribution of those paths. 
Moreover, as we will see later, this generalization relate to well-known sequences. So, a generalized Fine sequence on $[n]$ 
congruous $q$ modulo $p$ will be a Catalan path of length $2n$ such
that every primitive path ---subpaths that doesn't touch the 
$x$-axis--- starts with a rise at least equal to $q$ except the first one, who starts with $p$. Equivalently: Catalan paths such that 
the first rise is equal to $k\, q+p$, with $k<\lfloor \frac{n-p}{q} \rfloor$. We shall denote the set of generalized Fine sequences 
on $[n]$ congruous $q$ modulo $p$ by $F_n^{p,q}$.

As presented earlier in subsection \ref{sr}, we can also code the generalized Fine sequences with words $\omega$ on ${\bbbn}^n$.

\begin{definition}
\label{gf}
Generalized Fine sequences, $F_n^{p,q}$, are words $\omega = \omega_1
\omega_2 \ldots \omega_n$ on ${\bbbn}^n$, with $0 \leq p < q$, such that:
\begin{itemize}
\item[(i)] $\omega = \beta \gamma \, \; with \; \beta = 01 \ldots (p-1)$,
\item[(ii)] $\forall i \in [n-1], 0 \leq \omega_{i+1} \leq \omega_i + 1$,
\item[(iii)] $\forall i \in [p+1,n-q], \omega_i = 0 \; implies \; \omega_{i+1}\omega_{i+2} \ldots \omega_{i+q} = 12 \ldots (q-1)$.
\end{itemize}
\end{definition}

\begin{remark}
The following mappings give two bijections on generalized Fine
sequences: the first one between words and Catalan paths with
constraint on the primitive paths and the second one between the latter and Catalan
paths with constraint on the first rise. 
Let $x$ be a $(+1,+1)$ step and $\overline x$ a $(+1,-1)$ step.
\begin{itemize}
\item consider a word $\omega = \omega_1 \omega_2 \ldots \omega_n$ on
  ${\bbbn}^n$, satisfying Definition \ref{gf}. 
First, start the path with a $x$ step, then, for each $\omega_i$, if
  $\omega_i > \omega_{i-1}$ append a $x$ step, else, append
  $(\omega_{i-1} - \omega_i +1) \, \overline x$ steps and one $x$
  step. Finally, append $(\omega_n +1) \, \overline x$ steps. The reverse of this fonction is direct, thus giving us a bijection,
\item now, given a Catalan path with $k+1$ primitive paths: $x^p(\nu_p)\overline x \, x^q_{q_1}(\nu_{q_1}) \ldots
  x^q_{q_i}(\nu_{q_i}) \ldots x^q_{q_k} (\nu_{q_k})$, where
 $x^q_{q_i}(\nu_{q_i})$ is the $(i+1)^{th}$ primitive path and
  $\nu_{q_i}$ a postfix of this path. The corresponding path with
  constraint on the first rise will be: $x^{kq}\, x^p \, \overline x
  (\nu_p)\,(\nu_{q_1}) \ldots (\nu_{q_i}) \ldots (\nu_{q_k})$.
\end{itemize}
\end{remark}

\begin{figure}[h]
\center
\includegraphics{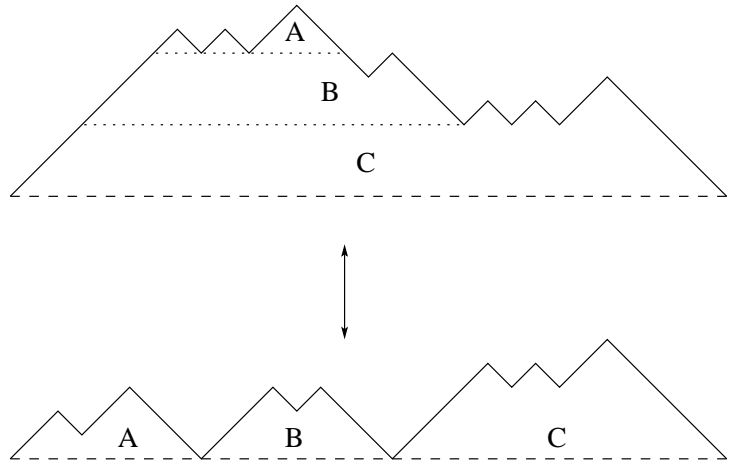}
\caption{Catalan paths corresponding to the generalized Fine sequence
  coded by $011201220123345$.}
\label{GFC}
\end{figure}

As an example, Figure \ref{GFC} illustrates the generalized Fine
sequence, with $q=3$ and $p=1$, coded by the integer sequence $011201220123345$.

The result which ensues from the above definitions is a direct enumerative formulae:

\begin{displaymath}
Card\{ F_n^{p,q} \} = \sum_{k=0}^{\lfloor \frac{n-p}{q} \rfloor}
\left\lbrack {2n-(kq+p)-1 \choose n-1} - {2n-(kq+p)-1 \choose n}
\right\rbrack
\end{displaymath}
where ${2n-k-1 \choose n-1} - {2n-k-1 \choose n}$ is the
$\alpha$-distribution of the Catalan numbers. That is a Catalan path
with a first rise of height $k$.

\section{Main theorems}
\label{mt}

We first state a general theorem about generalized Fine sequences. We then proceed with permutations with forbidden patterns.

\begin{definition}
\label{D1}
We consider the following succession system: \\

$
\left\{ \begin{array}{lll}
root = [P] \\
{[T]} & \longrightarrow  & [T],[3] \\
          & \underset{q-1} \longrightarrow & [q] \\
{[t]} & \longrightarrow & [T],[3],\ldots,[t+1] \\
\end{array} \right. \\
with \; P= \left\{ \begin{array}{ll}
p \; if \; p \geq 1\\
q \; else
\end{array} \right. \\
$

Take note that if $t < 2$, then the last succession rule only generates $[T]$. 
\end{definition}

\begin{theorem}
\label{tgfs}

Generalized Fine sequences, congruous p modulo q: $F_n^{p,q}$, can be charaterized by the succession system \ref{D1}, given just below.
\end{theorem}

\begin{theorem}
\label{der}

Derangements ---permutations such that $\forall i, \pi(i) \not = i$---
avoiding pattern $321$, which are in a one-to-one
correspondence with nonsingular similarity relations \cite{Deutsch},
can also be characterized by the succession system \ref{D1}, with
$p=0$ and $q=2$.
\end{theorem}

Among the generalized Fine sequences, two particular cases can be related to permutations with forbidden patterns. Those are:

\begin{theorem}
\label{fsf}

Nonsingular similarity relations, namely $F_n^{0,2}$, are in a
  one-to-one correspondence with three sets of permutations with
  forbidden patterns: $S_n({\cal F}_1)$, $S_n({\cal F}_2)$ and $S_n({\cal
  F}_3)$. With, 
${\cal F}_1 = \{1234,1243,1324,2134,2314,3124\}$,
${\cal F}_2 = \{1324,2134,2143,2314,3124,3214\}$ and \\
${\cal F}_3 =\{1342,2341,2413,2431,3142,3241\}$.

\end{theorem}

\begin{theorem}
\label{fsh}

Generalized Fine sequences congruous to \emph{one} modulo
\emph{three}, namely $F_n^{1,3}$, are in a one-to-one correspondence
with five sets of permutations with forbidden patterns: 
$S_{n-1}({\cal H}_1)$, $S_{n-1}({\cal H}_2)$, $S_{n-1}({\cal H}_3)$,
$S_{n-1}({\cal H}_4)$ and $S_{n-1}({\cal H}_5)$. With, ${\cal H}_1=\{1324,2314,2413,3124,3142,3214\}$, ${\cal
  H}_2=\{1234,1243,1324,1423,2314,3124\}$, ${\cal
  H}_3=\{2341,2413,2431,3412,3421,4231\}$, \\ 
${\cal H}_4=\{2134,2143,2314,3124,3214,4213\}$ and ${\cal H}_5
=\{1234,1243,1324,1423,2134,3124\}$.

\end{theorem}

\section{Some background}
\label{bk}

In order to follow the demonstration of the theorems, we give some
background and outline the sketch of the proves. But first of all,
let us define a usefull notation: say $_f \tau$ as ``the fordibben pattern
$\tau$''. This notation will come in handy thereafter. The reader already familiarized with the technique of
generating trees may now skip this section.

\subsection{Background on generating tree}

Originaly, the technique of \emph{generating tree} was introduced by F.R.K. Chung, R.L.~Graham, V.E.~Hoggat and M.~Kleiman
\cite{Chung} in order to enumerate Baxter permutations who avoided two patterns. Later, it was also applied to the study of various
permutations with forbidden subsequences by different authors (see for example 
\cite{Barcucci1,GireThese,GuibertThese,PergolaThese,WestPHD,WestCatalan,WestAG}).

A generating tree is a rooted, labeled tree such that the label of any
 vertex exclusively determined the number and the labels of his children.
Thus, any particular generating tree can be recursively defined by a
 \emph{succession system}, which is a set of succession rules
 consisting of a basis (the label of the root) and an inductive step
 (a set of labels), which are the children generated by any
 label. Moreover, each vertex sharing the same depth is associated
 with a combinatorial object of the same cardinality.

Consequently, if a succession system is shared by different generating
trees, the combinatorial objects related to those trees are in a
direct bijection. We say that they are $charaterized$ by the same succession system. 

Moreover, any succession system can be used to obtain recurrence
relations from which one may compute a closed form counting the objects themselves.

The relationship between structural properties of rules and 
the rationality, algebraicity or transcendence of the corresponding 
generating function has been investigated \cite{Banderier}. See also \cite{Bousquet}.
This technique also permits the random generation of the objects considered 
\cite{Barcucci2} (see also \cite[section~2.4]{GuibertThese}).

\subsection{Sketch of the proofs}

In order to prove the main theorems we will have to show, many times,
that a given set of permutations with forbidden subsequences can be characterized by a given succession system.

Given a set ${\cal E}$ of forbidden patterns and a permutation $\pi$
in $ S({\cal E})= \displaystyle{\bigcup_{n \geq 0}} S_n({\cal E})$. To prove that $S({\cal E})$ can 
be characterized by a succession system we have to:
\begin{itemize}
\item give a generating tree of $S({\cal E})$ characterized by the succession system,
\item associate each permutation belonging to $S({\cal E})$ with a label of
  the generating tree. Thus, the set of labels must form a partition
  of $S({\cal E})$. Furthermore, each labelized permutation
  belonging to $S({\cal E})$ must have a unique father in the generating tree,
\item prove that given a labeled permutation, all its children will
  have the same labels as given by the generating tree. Moreover, they
  must be unique and belongs to $S({\cal E})$. We also need
  to prove that the root has the right label.
\end{itemize}
Alltogether, these points define a bijection between a path in the generating tree and a permutation in $S({\cal E})$.

\begin{remark}
In the following, we will not explicitly give the generating trees
associated with the succession systems: they can be obtained
thoroughly with the definition of the labels and by inserting $n+1$ in
the active sites from left to right. The resulting permutations will
have, orderwise, the labels given by the succession rules. Thereafter,
we will not explicitly prove that the set of labels form a correct
partition, that all childrens are unique nor that the root has the
right label, as this can be done at first glance. Finally, since new
permutations are obtained by insertion, determining the unique father
is trivialy done by suppressing the greatest integer. Consequently
this point will not be developed either.
\end{remark}

\begin{definition}
Given a set ${\cal E}$ of forbidden patterns and $\pi$ a permutation
in $S_n({\cal E})$. We call a \emph{site} $i$ \emph{active} if the
permutation resulting from the insertion of $n+1$ between $\pi(i-1)$
and $\pi(i)$ is in $S_{n+1}({\cal E})$. A \emph{site} is called
\emph{inactive} otherwise. Moreover, a site is \emph{always active} if
it is active for any $\pi$ in $S_n({\cal E})$. In the subsequent
figures, a dot will denote a inactive site and a blank an active one.
\end{definition}
\begin{remark}
Consequently, to prove the inactivity of a site $i$, we have to find a
subsequence of $\pi$, containing $n+1$ in position $i$ orderisomorph
to some forbidden pattern $\tau$ in ${\cal E}$. Note that $n+1$
corresponds to the term of greatest ordinality in $\tau$.
\end{remark}
\begin{example}
Set $\tau=4312$ a forbidden pattern and $\pi = 326415$. The first
three sites are inactive as the subsequence 7615 is orderisomorph to
4312. Sites from four to seven are active as no subsequences in the
resulting permutations will be forbidden. Note that the last three
sites are always active as the value of greatest ordinality in $\tau$ is
in the first position.
\end{example}

\section{Generalized Fine sequences can
  be characterized by the succession system \ref{D1}}
\label{ptgfs}

\begin{definition}
\label{gtgf}
We consider the following generating tree:

$
\left\{ \begin{array}{lll}
root = 01\ldots (p-1)[P] \\
w=w'1[T] & \longrightarrow  & w1[T],w2[3] \\
          & \underset{q-1} \longrightarrow & w'01\ldots (q-1)[q] \\
w\not=w'1[t] & \longrightarrow & w1[T],w2[3],\ldots, wt[t+1] \\
\end{array} \right. \\
$
with $w$ labelized:
\begin{itemize}
\item $[T]$ if $w_n=1$,
\item $[t]$ otherwise.
\end{itemize}
\end{definition}

In order to prove Theorem \ref{tgfs}, we have to show that the following points hold true:
\begin{enumerate}
\item[$(i)$] exclusivity: considering an object generated by the tree in $F_n^{p,q}$, all of its
children belong to $ F_{m}^{p,q}$ with $m > n$. Moreover, the root is
in $F_p^{p,q}$,
\item[$(ii)$] completeness: every object in $F_n^{p,q}$ is generated
  by the generating tree given in Definition \ref{gtgf},
\item[$(iii)$] unicity: no object in $F_n^{p,q}$ appears more than once in the generating tree.
\end{enumerate}

\begin{proof}[Proof of Theorem \ref{tgfs}]
To prove the first point, we have to consider definitions \ref{gf} and
\ref{gtgf}. Clearly, all children of an object labeled $[t]$ 
fulfill condition $(ii)$ of Definition \ref{gf}. The same holds for
the objects labeled $[T]$, note that the second rule corresponds to
condition $(iii)$ of Definition \ref{gf}. Finally, condition $(i)$, the root belongs to
$F_p^{p,q}$. Now, given any generalized Fine sequences $\omega$ in
$F_n^{p,q}$, we can assign to it an unique labelized father with the
following mapping:
\begin{itemize}
\item $\omega_1\ldots\omega_{n-1}$ $[T]$ if ${\omega}_{n-1} = 1$,
\item $\omega_1\ldots\omega_{n-q}\,1$ $[T]$ if ${\omega}_{n-q+1}\ldots {\omega}_n =
  01\ldots(q-1)$,
\item $\omega_1 \ldots \omega_{n-1}$ $[\omega_n]$ otherwise.
\end{itemize}
This give a one-to-one correspondence between paths in the generating
tree and generalized Fine sequences and thus proves the second and third points.
\end{proof}

\section{Derangements avoiding $321$ are
  characterized by the succession system \ref{D1} with $p=1$ and $q=3$}
\label{sder}

Let's say $D_n(321)$ as $321$-avoiding derangements. 
\begin{definition}
\label{ssder}
We consider the following generating tree:

$
\left\{ \begin{array}{lll}
root = 21[2] \\
 {\pi(1)\ldots\pi(n-2)(n)\pi(n)[T]} 
& \rightarrow  & \pi(1)\ldots\pi(n-2)(n)(n+1)\pi(n)[T] \\
& \rightarrow  & \pi(1)\ldots\pi(n-2)\pi(n)(n+1)(n)[2] \\
& \rightarrow  & \pi(1)\ldots\pi(n-2)(n+1)\pi(n)(n)[3] \\
 {\pi(1)\ldots(n)\ldots\pi(n)[t]} & \rightarrow & \pi(1)\ldots(n)\ldots(n+1)\pi(n)[T] \\
& \rightarrow & \pi(1)\ldots(n)\ldots(n+1)\ldots\pi(n-1)\pi(n)[n+2-\pi^{-1}(n+1)] \\
& & with \; \pi^{-1}(n)+1 < \pi^{-1}(n+1) \leq n-1 \\
& \rightarrow & \pi(1)\ldots(n+1)\ldots\pi(n)(n)[t+1] \\
\end{array} \right. \\
$
with $\pi$ labelized:
\begin{itemize}
\item $[T]$ if $\pi(n-1)=n$, and $\pi(n) \not = n-1$,
\item $[t]$ else, with $t=n+1-\pi^{-1}(n)$.
\end{itemize}

\end{definition}

\begin{lemma}
\label{lder}
The generating tree \ref{ssder}, just given, generate $321$-avoiding derangements.
\end{lemma}

\begin{proof} 
Set $\pi$ in $D_n(321)$:
\begin{itemize}
\item if $\pi$ is labeled $[T]$ the succession system generate three
  new derangements:
\begin{itemize}
\item $n+1$ is inserted in position $n$, we clearly obtain a
  derangement and, by Definition \ref{ssder}, it has label $[T]$ and still avoid $_f 321$,
\item $\pi(n)$ substitute $n$ in $\pi$, $n+1$ substitute $\pi(n)$ and
  $n$ is placed at the end of the derangement. As $\pi(n) \not = n-1$,
  no fixed point can appear. Now, the new derangement has label $[2]$ as $\pi^{-1}(n+1)=n$ and
  $\pi(n+1)=n$. Moreover it clearly avoid $_f 321$,
\item $n+1$ substitute $n$ in $\pi$ and $n$ is placed at the end of
  the derangement. Again, no fixed point can appear. So, from Definition \ref{ssder} the new derangement
  has label $[3]$ and it clearly avoid $_f 321$.
\end{itemize}
\item if $\pi$ is labeled $[t]$, with $\pi(i)=n$, the succession system generate $t$ new derangements:
\begin{itemize}
\item $n+1$ is inserted in position $n$, since $\pi(n) \not =
  n$, no fixed point can appear and from Definition \ref{ssder} the new derangement will have label $[T]$,
\item $n+1$ is inserted in position $k$ in $[i+1\ldots n-1]$. First,
  take note that $\pi(i+1)\ldots\pi(n)$ is a strictly increasing
  subsequence. As a consequence, the derangement obtained by the insertion of $n+1$ in position
  $k$ is still $_f 321$-avoiding. Now, if a fixed point, $k+1 \leq \pi(j)
  \leq n$ appear, then the subsequence $\pi(j)\ldots \pi(n)$ should
  not have been strictly increasing in $\pi$ since $\pi(n) \not =
  n$. It follows a contradiction and consequently the new permutation
  is still a derangement. Definition \ref{ssder} account for the label of those new derangements,
\item $n+1$ substitute $n$ in $\pi$ and $n$ is placed at the end of
  the derangement. As above, no fixed point can appear and the new derangement is still
  $_f 321$-avoiding. Since $n+1$ is in position $i$ in the new
  derangement, it has label $[t+1]$.
\end{itemize}
\end{itemize}
This point implies the exclusivity. Now, set $i=\pi^{-1}(n)$, we can assign a unique
father to a given $\pi$ in $D_n(321)$ with the following mapping:

If $\pi$ has label:

\begin{itemize}
\item $[T]$ and $\left\{ \begin{array}{l} 
\pi(n-2)=n-1 \mapsto \pi(1)\ldots\pi(n-2)\pi(n)[T]\\
\pi(n-2)\not =n-1 \mapsto \pi(1)\ldots\pi(n-2)\pi(n)[n-\pi^{-1}(n-1)]\\
\end{array} \right. \\ $
\item $[2] \mapsto \pi(1)\ldots\pi(n-3)(n-1)\pi(n-2) [T]$
\item $[3]$ and $\left\{ \begin{array}{ll} 
\pi(n-1)\not =n-2,\; \pi(n)=n-1 & \mapsto
\pi(1)\ldots\pi(n-3)(n-1)\pi(n-1)[T]\\
\pi(n-1)\pi(n)=(n-2)(n-1) & \mapsto
\pi(1)\ldots\pi(i-1)(n-1)\pi(i+1)\ldots\pi(n-1)[2]\\
\end{array} \right. \\ $
\item $[t \geq 3]$ and $\pi(n) \not =n-1 \mapsto \pi(1)\ldots\pi(i-1)\pi(i+1)\ldots\pi(n)[n-\pi(n-1)]$
\item $[t \geq 4]$ and $\pi(1)\ldots(n)\ldots(n-1) \mapsto \pi(1)\ldots\pi(i-1)(n-1)\pi(i+1)\ldots\pi(n-1)[t-1]$
\end{itemize}
 
This mapping, along with the generating tree give us a one-to-one
correspondence between a path in the generating tree and a derangement in
$D_n(321)$. Since the root is trivialy in $D_1(321)$ the unicity and completness are achevied.
\end{proof}

\begin{proof}[Proof of Theorem \ref{der}]
 
First, recall that the succession system \ref{D1}, with $p=0$ and
$q=2$, is the following:

$
\left\{ \begin{array}{lll}
root = [2] \\
{[T]} & \rightarrow  & [T],[2],[3] \\
{[t]} & \rightarrow & [T],[3],\ldots,[t+1] \\
\end{array} \right. \\
$
As one can see, this succession system is isomorph to the succession
system charaterizing the generating tree \ref{ssder}. This remark and
Lemma \ref{lder} finish the
proof and give us a bijection between non singular similarity
relations and $321$-avoiding derangements.
\end{proof}

\section{Nonsingular similarity relations
  are in a one-to-one correspondence with three sets of forbidden patterns}
\label{pfsf}

\begin{lemma}
The generating trees of ${\cal F}_1, {\cal F}_2$ and ${\cal F}_3$ can
all be characterized by the succession system given by Definition
\ref{D1}, with $p=0$ and $q=2$.
\end{lemma}

\begin{proof}[Proof of Theorem \ref{fsf}]
The proof follows directly from the above
lemma. In order to prove the latter, we have to define generating trees
associated with ${\cal F}_1$, ${\cal F}_2$, ${\cal F}_3$ whose
succession systems are isomorph with the one given in Definition
\ref{D1}, with $p=0$ and $q=2$. This is done hereafter.
\end{proof}

\subsection{Generating tree of $S_n({\cal F}_1)$}

\begin{definition}
\label{ssf1}
\label{gtf1}
We consider the following succession system and labeling:

$
\left\{ \begin{array}{lll}
root = [2] \\
 {[T]} & \rightarrow  & [3],[T],[2] \\
 {[t]} & \rightarrow & [t+1],[T],[3],\ldots,[t] \\
\end{array} \right. \\
$

with $t$ an integer. \\
Given $\pi$ in $S_n({\cal F}_1)$, $\pi$ as label:
\begin{itemize}
\item $[T]$ if $\pi(2)=n$,
\item $[t]$ else, with $t$ the number of active sites of $\pi$.
\end{itemize}
\end{definition}

\begin{property}
\label{p1f1}

Given $\pi$ in $S_n({\cal F}_1)$:
\begin{itemize}
\item[(i)] take an active site, say $k$, if $k+1$ is inactive, then
  all sites from the first to $k$ are active and all others are
  inactive,
\item[(ii)] if $\pi$ has label $[T]$, then only the first three sites are active, 
\item[(iii)] if $\pi$ has label $[t]$ and $n+1$ is inserted in an
  active site $k$ greater than two, then, in the resulting
  permutation, site $k$ is active and $k+1$ inactive.
\end{itemize}
\end{property}

\begin{proof} $ $
\begin{itemize}
\item[(i)] results directly from the structure of the forbidden subsequences in ${\cal F}_1$,
\item[(ii)] set $\pi$ with label $[T]$ and consider the following
  permutation, $\pi(1)(n)\pi(3)(n+1)\ldots \pi(n)$, obtained by the
  insertion $n+1$ in its fourth site. Clearly, this permutation is
  forbidden as $\pi(1)(n)\pi(3)(n+1)$ will either be orderisomorph to
  $_f 2314$ or $_f 1324$. Now, the third site is active as neither
  $1342$ nor $2341$ are in ${\cal F}_1$. Applying Property \ref{p1f1}
  (i) finished the proof,
\item[(iii)] set $\pi$ with label $[t]$. The insertion of $n+2$ in
  position $k+1$ will give the subsequence $\pi(1)\pi(2)(n+1)(n+2)$
  orderisomorph to either $_f 2134$ or $_f 1234$. Now, site $k$
  remains active as if not, then some subsequence $a\,b\,(n+2)(n+1)$,
  with $a<b$, must be orderisomorph to $_f 1243$. Hence a
  contradiction as subsequence $a\,b\,(n+1)$ along with $n$ will have been
  orderisomorph either to $_f 3124$, $_f 1324$, $_f 1234$ or
  $_f1243$. And consequently $k$ shouldn't have been an active site in
  $\pi$.
\end{itemize} 
\end{proof}
Now, we prove that given a labelized permutation, its children will have the labels provided by the generating tree.

\begin{proof}
Set $\pi$ in $S_n({\cal F}_1)$:
\begin{itemize}
\item if $\pi$ is labeled $[T]$. As proven in Property \ref{p1f1} $(ii)$, $\pi$ has three active sites. Now, if $n+1$ is inserted in:
\begin{itemize}
\item the first site: the new permutation has label $[3]$. Indeed,
  insertion of $n+2$ in the fourth site will result in the subsequence
  $(n+1)\pi(1)(n)(n+2)$ being orderisomorph to $_f 3124$. However,
  site three is active as neither $3142$ nor $3241$ are forbidden
  patterns,
\item the second site: from Definition \ref{gtf1}, the new
  permutation has label $[T]$,
\item the third site: the new permutation has label $[2]$. Indeed,
  the first and second sites are active and the third is inactive as
  $\pi(1)(n)(n+2)(n+1)$ is orderisomorph to $_f 1243$.
\end{itemize}
\item if $\pi$ is labeled $[t]$ it has $t$ active sites. Now, if $n+1$ is inserted in:
\begin{itemize}
\item the first site: the new permutation has label $[t+1]$ since all
  sites of $\pi$ will remain active at which we add the new first
  site,
\item the second site: the new permutation has label $[T]$,
\item site $k$ in $[3,\ldots,t]$: using Property \ref{p1f1} $(iii)$ the new permutation will have label $[k]$.
\end{itemize}
\end{itemize}
\end{proof}

\subsection{Generating tree of $S_n({\cal F}_2)$}

\begin{definition}
\label{ssf2}
\label{gtf2}

We consider the following succession system and labeling:

$
\left\{ \begin{array}{lll}
root = [2] \\
 {[T]} & \rightarrow  & [T],[3],[2] \\
 {[t]} & \rightarrow & [T],[3],\ldots,[t+1] \\
\end{array} \right. \\
$

with $t$ an integer. \\
Given $\pi$ in $S_n({\cal F}_2)$, $\pi$ has label:
\begin{itemize}
\item $[T]$ if $\pi(1)=n$,
\item $[t]$ else, with $t$ the number of active sites of $\pi$.
\end{itemize}
\end{definition}

\begin{property}
\label{p1f2}
Given $\pi$ in $S_n({\cal F}_2)$:
\begin{itemize}
\item[(i)] take an active site, say $k$ ; if $k+1$ is inactive, then
  all sites from the first to $k$ are active and all others are
  inactive,
\item[(ii)] if $\pi$ has label $[T]$ ; then only the first three sites are active, 
\item[(iii)] if $\pi$ has label $[t]$ and $n+1$ is inserted in an
  active site $k$ greater than two ; then, in the resulting
  permutation, site $k+1$ is active and $k+2$ is inactive.
\end{itemize}
\end{property}

\begin{proof}
$ $
\begin{itemize}
\item[(i)] results directly from the structure of the forbidden subsequences in ${\cal F}_2$,
\item[(ii)] set $\pi$ with label $[T]$ and consider the following
  permutation: $(n)\pi(2)\pi(3)(n+1)\ldots \pi(n)$, obtained by the
  insertion $n+1$ in its fourth site. Clearly, this permutation is
  forbidden as $(n)\pi(2)\pi(3)(n+1)$ will either be orderisomorph to
  $_f 3124$ or $_f 3214$. Now, the third site is active as no patterns
  in ${\cal F}_2$ are orderisomorph to $1342$ or $2341$. Property \ref{p1f2} $(i)$, finishes the proof,
\item[(iii)] set $\pi$ with label $[t]$. The insertion of $n+2$ in
  site $k+2$ will give the subsequence $\pi(1)(n+1)\pi(k)(n+2)$
  orderisomorph to either $_f 1324$ or $_f 2314$. Now, site $k+1$ is
  active as if not, then some subsequence $a\, b\, (n+1)(n+2)$, with
  $a>b$, must be orderisomorph to $_f 2134$. Hence a contradiction
  appears, as subsequence $a\,b\,(n+1)$ along with $n$ will have been
  orderisomorph either to $_f 2143$, $_f 2134$, $_f 2314$ or $_f
  3214$ and consequently $k$ should not have been an active site in $\pi$.
\end{itemize}
\end{proof}
Now, we prove that given a labelized permutation, its children will have the labels provided by the generating tree.

\begin{proof}
Set $\pi$ in $S_n({\cal F}_2)$:
\begin{itemize}
\item if $\pi$ is labeled $[T]$. As proven in Property \ref{p1f2} $(ii)$, $\pi$ has three active sites. Now, if $n+1$ is inserted in:
\begin{itemize}
\item the first site: the new permutation has label $[T]$,
\item the second site: then the new permutation has label
  $[3]$. Indeed, the fourth site is inactive as the subsequence
  $(n)(n+1)\pi(2)(n+2)$ will be orderisomorph to $_f 2314$. Moreover, the third site is active as pattern $2341$ is allowed,
\item the third site: the new permutation has label $[2]$. Indeed,
  the first and second sites are active and the third is inactive as $(n)\pi(2)(n+2)(n+1)$ is orderisomorph to $_f 2143$.
\end{itemize}
\item if $\pi$ is labeled $[t]$ it as $t$ active sites. Now, if $n+1$ is inserted in:
\begin{itemize}
\item the first site: the new permutation as label $[T]$,
\item site $k$ in $[2,\ldots,t]$. Using Property \ref{p1f2} $(iii)$ we see that the new permutation will have label $[k+1]$.
\end{itemize}
\end{itemize}
\end{proof}

\subsection{Generating tree of $S_n({\cal F}_3)$}

\begin{definition}
\label{ssf3}
\label{gtf3}
We consider the following succession system and labeling:

$
\left\{ \begin{array}{lll}
root = [2] \\
 {[T]} & \rightarrow  & [T],[2],[3] \\
 {[t]} & \rightarrow & [T],[3],\ldots,[t+1] \\
\end{array} \right. \\
$

with $t$ an integer. \\ 
Given $\pi$ in $S_n({\cal F}_3)$, $\pi$ has label:
\begin{itemize}
\item $[T]$ if $\pi(1)=n$,
\item $[t]$ else, with $t$ the number of active sites of $\pi$.
\end{itemize}
\end{definition}

\begin{figure}[h]
\center
\input{Schema_F3}
\caption{generating tree of $S_n({\cal F}_3)$.}
\label{F3}
\end{figure}
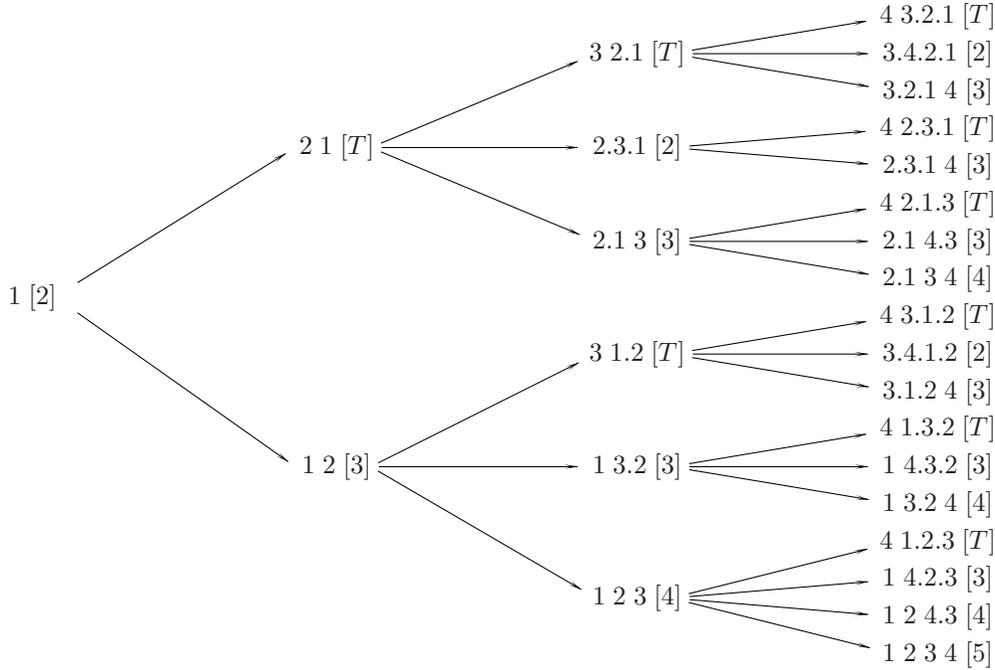

\begin{property}
\label{p1f3}

Given $\pi$ in $S_n({\cal F}_3)$:
\begin{itemize}
\item[(i)] The first and last sites are always active,
\item[(ii)] if $\pi$ has label $[T]$, then only the first, the second and last sites are active,
\item[(iii)] if $\pi$ has label $[t]$ and $n+1$ is inserted in site
  $1<k<n+1$, then all sites from $2$ to $k$ remain active and all sites from $k+1$ to $n$ are inactive.
\end{itemize}
\end{property}

\begin{proof} $ $
\begin{itemize}
\item[(i)] results directly from the structure of the forbidden subsequences in ${\cal F}_3$,
\item[(ii)] the second site is active as, with $\pi(1)=n$, patterns
  $_f 2413$ and $_f 2431$ cannot appear. Moreover, all sites from
  three to $n$ are inactive as subsequence $(n)\pi(2)(n+1)\pi(n)$ will be either orderisomorph to $_f 3142$ or $_f 3241$,
\item[(iii)] consider the following permutation,
  $\pi(1)\ldots(n+1)\dots\pi(n)$, obtained after insertion of $n+1$ in
  $\pi$ in position $k$. All sites from $k+1$ and $n$ are inactive as
  the subsequence $\pi(1)(n+1)(n+2)\pi(n)$ is orderisomorph either to
  patterns $_f 1342$ or $_f 2341$. Now, if a site $l$ before $k$
  becomes inactive, it must contain a subsequence orderisomorph to
  patterns $_f 2413$ or $_f 2431$, with $n+1$ and $n+2$ respectively
  in position $k$ and $l$. But then, the subsequence $241$ along with
  $n$ will have been orderisomorph to one of the following patterns:
  $_f 3241$, $_f 2341$, $_f 2431$ or $_f 2413$. Hence a
  contradiction. So $l$ remains active.
\end{itemize}
\end{proof}
Now, we prove that given a labelized permutation, its children will have the labels provided by the generating tree.

\begin{proof}
Set $\pi$ in $S_n({\cal F}_3)$:
\begin{itemize}
\item if $\pi$ is labeled $[T]$. As proven in Property \ref{p1f3} $(ii)$, $\pi$ has three active sites. Now, if $n+1$ is inserted in:
\begin{itemize}
\item the first site: the new permutation has label $[T]$,
\item the second site: the new permutation has label $[2]$. Indeed,
  the first and last sites are always active, as stated in Property
  \ref{p1f3} $(i)$, moreover, site two is inactive as
  $(n)(n+2)(n+1)\pi(n)$ is orderisomorph to $_f 2431$. Finally, all
  sites from three to $n$ are inactive as $(n)(n+1)(n+2)\pi(n)$ is
  orderisomorph to $_f 2341$.
\end{itemize}
\item if $\pi$ is labeled $[t]$, it has at least two active sites. Now, if $n+1$ is inserted in:
\begin{itemize}
\item the first site: the new permutation has label $[T]$,
\item the $k^{th}$ active site, with $k>1$ , then Properties \ref{p1f3}
  $(ii)$ and $(iii)$ imply that the new permutation has label $[k+1]$.
\end{itemize}
\end{itemize}
\end{proof}

\newpage

\section{Generalized Fine sequences
  congruous three modulo one are in one-to-one correspondence with
  five sets of forbidden patterns}
\label{pfsh}

We begin by stating three lemmas which are proven in the next sections.

\begin{lemma}
\label{lss1}
Generalized Fine sequences congruous three modulo one, namely $F_n^{1,3}$, are in a
one-to-one correspondence with permutations with forbidden
patterns $S_{n-1}({\cal H}_1)$.
\end{lemma}

\begin{lemma}
\label{lss2}
Permutations with forbidden patterns $S_n({\cal E})$, with ${\cal
  E}$ taking value in $\{{\cal H}_1^{\star}$, ${\cal H}_2$, ${\cal
  H}_3$, ${\cal H}_4 \}$  can be characterized by the same succession system, thus giving a bijective correspondence among them.
\end{lemma}

\begin{lemma}
\label{lss3}
Permutations with forbidden patterns $S_n({\cal H}_3^{-1\,c})$ and
$S_n({\cal H}_5)$ can be characterized by the same succession system,
thus giving a bijective correspondence between them.
\end{lemma}

\begin{figure}[h]
\center
\input{Schema_preuve}
\caption{an overview of the correspondences between the five sets.}
\label{corr}
\end{figure}

\begin{proof}[Proof of Theorem \ref{fsh}]
Lemma \ref{lss1} establish a bijection between
$F_n^{1,3}$ and $S_{n-1}({\cal H}_1)$. Now, the operation of $mirror$
---denoted by $\star: \forall i \in n,\, \pi^{\star}(i)=\pi(n-i+1)$--- on ${\cal H}_1$
 gives also a one-to-one correspondence between
$S_n({\cal H}_1)$ and $S_n({\cal H}_1^{\star})$. It follows that all
permutations with forbidden patterns given in Lemma \ref{lss2}
are, through $S_{n-1}({\cal H}_1^{\star})$ and $S_{n-1}({\cal H}_1)$, in a
one-to-one correspondence with $F_n^{1,3}$. Finally, the last set,
${\cal H}_5$, is obtained  with the operations of $inverse$ ---denoted by
$-1$--- and $complementary$ ---denoted
by $c:\forall i \in n,\,\pi^c(i)-\pi(i)=n+1$--- on ${\cal H}_3$ along with Lemma \ref{lss3}.
\end{proof}

\subsection{A bijection
  between $S_{n-1}({\cal H}_1)$ and $F_n^{1,3}$}

\subsubsection{Generating tree of $S_n({\cal H}_1)$}

\begin{definition}
\label{gth1}
We consider the following generating tree:

$
\left\{ \begin{array}{lll}
root = 1[T] \\
 {(n)\pi(2)\ldots\pi(n)[T]} & \rightarrow  & (n+1)(n)\pi(2)\ldots\pi(n)[T],(n)(n+1)\pi(2)\ldots\pi(n)[3] \\
       & \underset{2} \rightarrow & (n+1)(n)(n+2)\pi(2)\ldots\pi(n)[3] \\
 { \pi(1)\ldots\pi(n)[t]} & \rightarrow & (n+1)\pi(1)\ldots\pi(n)[T] \\
       & \rightarrow & \pi(1)(n+1)\ldots\pi(n)[3],\ldots,\pi(1)\ldots\pi(\pi^{-1}(n))(n+1)\ldots\pi(n)[t+1] \\
\end{array} \right. \\
$
with $\pi$ labelized:
\begin{itemize}
\item $[T]$ if $\pi(1)=n$,
\item $[t]$ else, with $t$ the number of active sites of $\pi$.
\end{itemize}
\end{definition}

\begin{property}
\label{p1h1}
Given $\pi$ in $S_n({\cal H}_1)$:
\begin{itemize}
\item[(i)] if $\pi$ has label $[T]$ and $\pi(2) \not =n-1$ then only
  the first two sites are active. Moreover, if $\pi(2)=n-1$ then only the first three sites are active,
\item[(ii)] if $\pi$ has label $[t]$. Set $i$ such that
  $\pi(i)=n$. Then all sites from $1$ to $i+1$ are active except the
  second one if $\pi(1)\pi(2)\pi(3)$ is orderisomorph to
  $213$. Moreover, $t>2$.
\end{itemize}
\end{property}

\begin{proof} $ $
\begin{itemize}
\item[(i)] results directly from the structure of the forbidden subsequences in ${\cal H}_1$,
\item[(ii)] the first site is trivialy active as no pattern in ${\cal
  H}_1$ begins with 4. Now, if a site $j$ in $[3\ldots i]$ is inactive,
  then some subsequence $\pi(k)(n+1)\pi(l)(n)$, with $1<k<l<i+1$ must
  be orderisomorph to the forbidden pattern $_f 2413$. Then, either
  $\pi(1)>\pi(k)$ and $\pi(1)\pi(k)\pi(l)\pi(i)$ is orderisomorph to
  $_f 3214$, either $\pi(1)<\pi(k)$ and the same subsequence is now
  orderisomorph to $_f 2314$. Hence a contradiction since those two
  patterns are forbidden. Site $i+1$ is active as no patterns in
  ${\cal H}_1$ end with 34. All sites $j$ in $[i+2\ldots n+1]$ are
  inactive as patterns $_f 1324$ and $_f 2314$ are forbidden. It
  follows, since $i$ is greater than one, that the first, third and
  fourth sites are always active and, consequently, $t>2$. Finally, as
  $_f 2413$ is a forbidden pattern, the second site is inactive if
  $\pi(1)\pi(2)\pi(3)$ is orderisomorph to $213$.
\end{itemize}
\end{proof}

\begin{lemma}
\label{lh1}
The generating tree \ref{gth1}, just given, generate $S_n({\cal H}_1)$.
\end{lemma}

\begin{proof} First, we prove that given a labelized permutation, its children will
have the labels given by the generating tree, this point implies exclusivity.

Set $\pi$ in $S_n({\cal H}_1)$:
\begin{itemize}
\item if $\pi$ is labeled $[T]$. As proven in Property \ref{p1h1}
  $(i)$, $\pi$ has always two active sites. Three succession rules
  apply:
\begin{itemize}
\item $n+1$ is inserted in the first site: from Definition \ref{gth1}, the resulting
  permutation has label $[T]$,
\item $n+1$ is inserted in the second site: the resulting permutation corresponds to label
  $[3]$ as all sites $i$ in $[4\ldots n+2]$ are inactive since
  $(n)(n+1)\pi(2)(n+2)$ is orderisomorph to $_f 2314$ and
  $(n)(n+1)\pi(2)$ isn't orderisomorph to $213$ (Property \ref{p1h1} $(ii)$).
\item $n+1$ is inserted in the first site and $n+2$ in the second one: from Property \ref{p1h1} $(ii)$, the resulting permutation,
  $(n+1)(n)(n+2)\pi(2)\ldots\pi(n)$, has label $[3]$.
\end{itemize}
\item if $\pi$ is labeled $[t]$, it has $t$ active sites. Now, if
  $n+1$ is inserted in:
\begin{itemize}
\item the first site: from Definition \ref{gth1}, the resulting permutation has label $[T]$,
\item site $k$ in $[3\ldots t+1]$: using Property \ref{p1h1} $(ii)$, we see that the resulting permutations will have labels
  ranging from $[3]$ to $[t+1]$.
\end{itemize}
\end{itemize}
Now, we can assign a unique father to a given $\pi$ in $S_n({\cal
  H}_1)$ with the following mapping: \\
\medskip
$
\left\{ \begin{array}{l}
\pi(2)\pi(4)\ldots\pi(n)[T]$ if $\pi(1)=n-1$ and $\pi(3)=n \\
\pi(1)\ldots\pi(\pi^{-1}(n)-1)\pi(\pi^{-1}(n)+1)\ldots\pi(n)$
  else, with the label given by Definition \ref{gth1}.$
\end{array} \right. \\
$
\medskip
This point give us both unicity and compleness.
\end{proof}

\begin{proof}[Proof of Lemma \ref{lss1}] First, remark that the succession system associated to the generating
tree \ref{gth1} is isomorph to the succession system \ref{D1}
associated with $F_n^{1,3}$ with one application of the succession
rules on the root. The last point account for
the cardinality's shift between the objects. Finally, Lemma \ref{lh1}
finish the proof.
\end{proof}

\subsection{Permutations avoiding ${\cal H}_1^{\star}$, ${\cal
  H}_2$, ${\cal H}_3$ and ${\cal H}_4$ are all characterized by the same succession system}

\begin{proof}[Proof of Lemma \ref{lss2}] All the succession systems
  associated with the generating trees defined hereafter are
  equivalent up to an isomorphism. As those generating trees produces $S_n({\cal H}_1^{\star})$, $S_n({\cal
  H}_2)$, $S_n({\cal H}_3)$ and $S_n({\cal H}_4)$, they are all
  characterized by the same succession system.
\end{proof}

\subsubsection{Generating tree of $S_n({\cal H}_1^{\star})$}
Remember that ${\cal H}_1^{\star}$ is the $mirror$ of ${\cal H}_1$. So we have ${\cal H}_1^{\star}=\{2341,2413,2431,4231,3142,3241\}$.

\begin{definition}
\label{ssh1*}
\label{dseH1*}
We consider the following succession system and labeling:

$
\left\{ \begin{array}{lll}
root = [A,2] \\
{[A,t]} & \rightarrow & [B,3],\ldots,[B,t+1],[A,t+1] \\
{[B,t]} & \rightarrow & [B,3],\ldots,[B,t],[A,1],[A,t] \\
\end{array} \right. \\
$
Given $\pi$ in $S_n({\cal H}_1^{\star})$ and $i,j$ such that
$\pi(i)=n-1$ and $\pi(j)=n$, $\pi$ has label:
\begin{itemize}
\item $[A,t]$ if $i<j$ and $\pi$ has $t$ active sites,
\item $[B,t]$ if $j<i$ and $\pi$ has $t$ active sites.
\end{itemize}
\end{definition}

\begin{figure}[h]
\center
\input{Schema_H1star}
\caption{generating tree of $S_n({\cal H}_1^{\star})$.}
\label{H1star}
\end{figure}

\begin{property}
\label{p1h1*}
Set $\pi$ in $S_n({\cal H}_1^{\star})$:
\begin{itemize}
\item[(i)] if $\pi$ has label $[A,t]$ and $j<n$, then the only active site is $n+1$,
\item[(ii)] if $\pi$ has label $[A,t]$ and $j=n$, then sites $n$,
  $n+1$ are always active, all sites from $j+1$ to $n-1$ are
  inactive. Moreover if $n+1$ is inserted in position $n+1$, all
  active sites in $\pi$ remain active in the resulting permutation,
\item[(iii)] if $\pi$ has label $[B,t]$. Then all sites from $j+2$ to
  $n$ are inactive. Sites $j$, $j+1$ and $n+1$ are always
  active. Moreover if $n+1$ is inserted in a position $k$ less than
  $j+1$, all active sites belonging to $[1\ldots k]$ in $\pi$ remain
  active in the resulting permutation.
\end{itemize} 
\end{property}

\begin{proof} $ $
\begin{itemize}
\item[(i)] By hypothesis, $\pi = \pi(1)\ldots \pi(i) \ldots \pi(j)
  \ldots \pi(n)$. Now, insertion of $n+1$ in sites $[1\ldots i]$,
  $[i+1\ldots j]$ and $[j+1\ldots n]$ will result in subsequences
  being respectively orderisomorph to the forbidden patterns $_f
  4231$, $_f 2431$ and $_f 2341$. Therefore, as no forbidden patterns
  in ${\cal H}_1^{\star}$ end with 4, site $n+1$ is active,
\item[(ii)] sites $n$, $n+1$ are always active as forbidden patterns
  end neither with $4$ nor $43$. Moreover the resulting permutations
  will keep their active sites as none of patterns $1423$, $4123$,
  $1432$ or $4132$ are forbidden. Now, consider site $k$ in
  $[j+1\ldots n-1]$, the permutations resulting from the insertion of
  $n+1$ in site $k$ will be forbidden as the subsequence
  $(n-1)(n+1)\pi(n-1)(n)$ is orderisomorph to $_f 2413$, a forbidden
  pattern,
\item [(iii)] consider $\pi = \pi(1)\ldots (n) \ldots (n-1) \ldots
  \pi(n)$. As patterns $_f 2413$ and $_f 4231$ must be avoided, it
  follows that
  $\pi(1)\ldots\pi(j-1)<\pi(j+1)\ldots\pi(i-1)<\pi(i+1)\ldots\pi(n)$.
  Consequently, all sites $k$ in $[j+2\ldots n]$ are inactive as
  $(n)\pi(k-1)(n+1)\pi(n)$ is orderisomorph either to $_f 3142$ or $_f
  3241$. Site $n+1$ is active as no forbidden patterns end with
  4. Now, sites $j$ and $j+1$ are active as neither pattern $_f 2341$
  nor $_f 2431$ can appear. Next, if an active site becomes inactive
  after insertion of $n+1$, then some subsequence containing $n$ and
  $n+1$ must be orderisomorph to the forbidden patterns $4231$, $2431$
  or $2413$. The first two cases are dismissed as no value before $n$
  is greater than after (see above). Next consider subsequence
  $\pi(1)(n+2)\pi(2)(n+1)$ orderisomorph to $_f 2413$. As
  $\pi(1)(n+1)\pi(2)(n)$ is orderisomorph to the forbidden pattern $_f
  2413$, this site would not have been active in $\pi$, hence a
  contradiction.
\end{itemize} 
\end{proof}
Now, we prove that given a labelized permutation, its children will have the labels provided by the generating tree.
\begin{proof}
Set $\pi$ in $S_n({\cal H}_1^{\star})$:
\begin{itemize}
\item if $\pi$ is labeled $[A,t]$, it has $t$ active sites. Now, if
  $n+1$ is inserted in:
\begin{itemize}
\item the last active site: if $j<n$, as stated in Property
  \ref{p1h1*} $(i)$, $t=1$. The resulting permutation,
  $\pi(1)\ldots(n-1)\ldots(n)\ldots\pi(n)(n+1)$, has label
  $[A,2]$. Indeed, the insertion of $n+2$ in sites, $[1\ldots i] \;,
      [i+1\ldots j]$ and $[j+1\ldots n]$, will, respectively, create
      subsequences orderisomorph to $_f 4231$, $_f 2431$ and $_f
      2413$. Now, if $j=n$, the resulting permutation, as a direct
      consequence of the structure of the forbidden patterns, has
      label $[A,t+1]$,
\item active sites from the first to the last but one, which imply
  $j=n$. From Property \ref{p1h1*} $(ii)$, and Definition
  \ref{dseH1*}, all permutations resulting from insertion of $n+1$ in
  all but the last active sites will be labelized $[B,k]$, with $k$
  taking values in $[3\ldots t+1]$. Indeed, all sites before the
  insertion are conserved and both sites, $j+1$ and $n+1$, are always
  active (Property \ref{p1h1*} $(iii)$).
\end{itemize}
\item if $\pi$ is labeled $[B,t]$, it has $t$ active sites. Now, if
  $n+1$ is inserted in:
\begin{itemize}
\item the first $(t-2)^{th}$ sites: again, we apply Property
  \ref{p1h1*} $(iii)$, and, consequently, those permutations will be
  labelized $[B,k]$, with $k$ taking values in $[3\ldots t+1]$,
\item site $j+1$: from Property \ref{p1h1*} $(i)$ the resulting
  permutation will have label $[A,1]$,
\item site $n+1$: the resulting permutation will have label $[A,t]$, since, from Property \ref{p1h1*}
  $(iii)$, the $t-2$ active sites before $n$ will be conserved and the
  two last sites are always active (Property \ref{p1h1*} $(ii)$).
\end{itemize}
\end{itemize}
\end{proof}

\subsubsection{Generating tree of $S_n({\cal H}_2)$}

\begin{definition}
\label{ssh2}
\label{dseH2}
We consider the following succession system and labeling:

$
\left\{ \begin{array}{lll}
root = [A,2] \\
{[A,t]} & \rightarrow & [A,t+1],[B,3],\ldots,[B,t+1] \\
{[B,t]} & \rightarrow & [A,t],[B,3],\ldots,[B,t],[A,1] \\
\end{array} \right. \\
$
Given $\pi$ in $S_n({\cal H}_2)$ and $i,j$ such that $\pi(i)=n-1$ and
$\pi(j)=n$, $\pi$ has label:
\begin{itemize}
\item $[A,t]$ if either $j=1$ or $1<i<j$ and $\pi$ has $t$ active sites,
\item $[B,t]$ if $1<j<i$ and $\pi$ has $t$ active sites.
\end{itemize}
\end{definition}

\begin{property}
\label{p1h2}
Set $\pi$ in $S_n({\cal H}_2)$:
\begin{itemize}
\item[(i)] if $\pi$ has $t$ active sites then, all sites from $t+1$ to
  $n+1$ are inactive and consequently, all sites in $[1\ldots t]$ are
  active. Moreover, if $\pi$ has label $[A,t]$ and $1<i<j$ then only
  the first site is active,
\item[(ii)] $\pi$ has label $[B,t]$ if and only if $\pi(t-1)=n$.
\end{itemize}
\end{property}

\begin{proof} $ $
\begin{itemize}
\item[(i)] if a site $k$ is inactive because some subsequence
  $a\,(n+1)\,b\,c$ is orderisomorph to the forbidden pattern $_f 1423$
  then all sites from $k+1$ to $n+1$ are inactive as $_f 1243$ and $_f
  1234$ are also forbidden. The same holds true if some subsequence
  $a\, b\, (n+1)\,c$ is orderisomorph to the forbidden pattern $_f
  1243$,
\item[(ii)] site $j+2$ is inactive as permutation
  $\pi(1)\pi(j)\pi(j+1)(n+1)$ is orderisomorph to the forbidden
  pattern 1324. Now, site $j+1$ must be active. If not, some
  subsequence, $a\, b\, c\, (n+1)$, $a\, b\, (n+1)\,c$ or $a\, (n+1)\,b
  \, c$ of the permutation obtained by insertion of $n+1$ in position
  $j+1$, shall be orderisomorph to, at least, one of the forbidden
  patterns. Then, remembering that $\pi(i)=n-1$ and $\pi(j)=n$, $a\,
  b\, c\,\pi(j)$, $a\, b\, \pi(j) \pi(i)$, $a\,b\,\pi(j)\, c$ or
  $a\,\pi(j) \,b\,c$ will have been, also, orderisomorph to the same
  patterns. As a consequence $\pi$ was not in $S_n({\cal H}_2)$, hence
  a contradiction. Next, using Property \ref{p1h2} $(i)$, all sites
  from 1 to $j+1$ are active. Finally, Definition \ref{dseH2} implies $\pi(t-1)=n$.
\end{itemize}
\end{proof}
Now, we prove that given a labelized permutation, its children will have the labels provided by the generating tree.
\begin{proof}
Set $\pi$ in $S_n({\cal H}_2)$:
\begin{itemize}
\item if $\pi$ is labeled $[A,t]$, it has $t$ active sites. Now, if
  $n+1$ is inserted in:
\begin{itemize}

\item the first active site $j \not =1$: then, only the first site is
  active as any other site will result in subsequences,
  $\pi(1)(n+1)\pi(i)\pi(j)$, $\pi(1)\pi(i)(n+1)\pi(j)$ and
  $\pi(1)\pi(i)\pi(j)(n+1)$ being all, respectively orderisomorph to
  $_f 1423$, $_f 1243$ and $_f 1234$. The resulting permutation has
  label $[A,2]$. Indeed, the second site is active as forbidden
  patterns begin with neither $342$ nor $341$ and the third site is
  inactive as $\pi(1)(n+2)\pi(i)\pi(j)$ is orderisomorph to $_f
  1423$. Now, if $j=1$, the resulting permutation has label
  $[A,t+1]$. Indeed, if site $t$ of $\pi$ becomes inactive then some
  subsequence $(n+1)\, a\, b\, (n+2)$, with $n+2$ in position $t+1$,
  must be orderisomorph to $_f 3124$. But then $\pi(j)\, a\, b\,
  (n+1)$, with $n+1$ in position $t$, will also, have been
  orderisomorph to the same pattern. Hence a contradiction. It follows
  that all active sites of $\pi$ remain active, to which, we add the
  new first site,
\item the second to last active sites: the resulting permutations
  have labels in $[B,3]\ldots [B,t]$ as a direct consequence of
  Definition \ref{dseH2} and Property \ref{p1h2} $(iii)$.
\end{itemize}
\item if $\pi$ is labeled $[B,t]$, it has $t$ active sites. Now, if $n+1$ is inserted in:
\begin{itemize}
\item the first active site: the resulting permutation has label
  $[A,t]$. Indeed, site $t$ becomes inactive as $(n+1)\pi(1)(n)(n+2)$
  is orderisomorph to $_f 3124$. Site $t-1$ remains active as, in
  contradiction, if some subsequence $(n+1)\, a\, b\, (n+2)$, with
  $n+2$ in the former position $t-1$, is orderisomorph to $_f 3124$
  then $a\, b\, \,\pi(j)\pi(i)$ would have been orderisomorph to $_f
  1243$,
\item the second to the $t-1^th$ active sites: this will result with
  permutations labelized from $[B,3]$ to $[B,t]$ as a direct
  consequence of Definition \ref{dseH2} and Property \ref{p1h2}
  $(iii)$,
\item the last active site: the resulting permutation,
  $\pi(1)\ldots(n)(n+1)\ldots\pi(n)$, has label $[A,1]$ as seen
  previously.
\end{itemize}
\end{itemize}
\end{proof}

\subsubsection{Generating tree of $S_n({\cal H}_3)$}

\begin{definition}
\label{ssh3}
\label{dseH3}
We consider the following succession system and labeling:

$
\left\{ \begin{array}{lll}
root = [A,2] \\
{[A,t]} & \rightarrow & [B,3],\ldots,[B,t+1],[A,t+1] \\
{[B,t]} & \rightarrow & [B,3],\ldots,[B,t],[A,1],[A,t] \\
\end{array} \right. \\
$
Given $\pi$ in $S_n({\cal H}_3)$ and $i,j$ such that $\pi(i)=n-1$ and
$\pi(j)=n$, $\pi$ has label:
\begin{itemize}
\item $[A,t]$ if $i>j$ and $\pi$ has $t$ active sites,
\item $[B,t]$ if $i<j$ and $\pi$ has $t$ active sites.
\end{itemize}
\end{definition}

\begin{property}
\label{p1h3}
Set $\pi$ in $S_n({\cal H}_3)$:
\begin{itemize}
\item[(i)]  if $\pi$ has label $[A,t]$ and $j<n$, then the only active site is $n+1$,
\item[(ii)] if $\pi$ has label $[A,t]$ and $j=n$, then sites $n$,
  $n+1$ are always active, all sites from $j+1$ to $n-1$ are
  inactive. Moreover if $n+1$ is inserted in position $n+1$, all
  active sites in $\pi$ remain active in the resulting permutation,
\item[(iii)] if $\pi$ has label $[B,t]$. Then all sites from $j+1$ to
  $n-1$ are inactive. Sites $j$, $n$ and $n+1$ are always
  active. Moreover if $n+1$ is inserted in a position $k$ less than
  $j+1$, all active sites belonging to $[1\ldots k]$ in $\pi$ remain
  active in the resulting permutations.
\end{itemize} 
\end{property}

\begin{proof} $ $
\begin{itemize}
\item[(i)] by hypothesis, $\pi = \pi(1)\ldots \pi(i) \ldots \pi(j)
  \ldots \pi(n)$. Now, insertion of $n+1$ in sites $[1\ldots i]$,
  $[i+1\ldots j]$ and $[j+1\ldots n]$ will result in subsequences
  being respectively orderisomorph to $_f 4231$, $_f 2431$ and $_f
  2341$. Therefore, as no forbidden patterns in ${\cal H}_3$ end with
  $4$, site $n+1$ is active,
\item[(ii)] sites $n$, $n+1$ are always active as forbidden patterns
  end neither with $4$ nor $43$. Moreover the resulting permutations
  will keep its active sites as patterns $1423$, $4123$, $1432$, $4132$
  are not forbidden. Now, consider site $k$ in $[j+1\ldots n-1]$, the
  permutation resulting from the insertion of $n+1$ in site $k$ will
  be forbidden as the subsequence $(n-1)(n+1)\pi(n-1)(n)$ is
  orderisomorph to $_f 2413$,
\item[(iii)] consider $\pi = \pi(1)\ldots (n) \ldots (n-1) \ldots
  \pi(n)$. As patterns $_f 2413$ and $_f 4231$ must be avoided it
  follows that
  $\pi(1)\ldots\pi(j-1)<\pi(j+1)\ldots\pi(i-1)<\pi(i+1)\ldots\pi(n)$.
  Consequently, all sites $k$ in $[j+1\ldots n-1]$ are inactive as
  $(n)(n+1)\pi(n-1)\pi(n)$ is orderisomorph either to $_f 3412$ or $_f
  3421$. Site $n+1$ is active as no forbidden patterns end with
  4. Site $n-1$ is active, as supposing it generates a subsequence
  $b\,c\,(n+1)\,a$ orderisomorph to $_f 2341$ ; then the same
  subsequence with $\pi(j)$ would have been orderisomorph to $_f
  4231$, $_f 2431$ or $_f 2341$. Site $j$ is active as pattern $_f
  2431$ cannot appear. Now, if $n+1$ is inserted in position $k<j+1$,
  then all active sites before $k$ remain active. In contrary, if such
  a site becomes inactive, then some subsequence, involving $n+1$ and
  $n+2$ must be orderisomorph to $_f 2413$, $_f 2431$ or $_f
  4231$. For the two latter cases we have a direct contradiction as
  all values, except $n+2$, before $n+1$ are lesser than those after
  $n+1$. For the former, the site should not have been active as the
  same subsequence, along with $n+1$ in position $k$ and $n$ would
  have been orderisomorph to $_f 2413$.
\end{itemize} 
\end{proof}
Now, we prove that given a labelized permutation, its children will have the labels provided by the generating tree.
\begin{proof}
Set $\pi$ in $S_n({\cal H}_3)$:

\begin{itemize}
\item if $\pi$ is labeled $[A,t]$, it has $t$ active sites. Now, if
  $n+1$ is inserted in:
\begin{itemize}
\item the first $(t-1)^th$ active sites: if $j<n$, as stated in
  Property \ref{p1h3} $(ii)$, $t=1$. The resulting permutation,
  $\pi(1)\ldots(n-1)\ldots(n)\ldots\pi(n)(n+1)$, has label
  $[A,2]$. Indeed, the insertion of $n+2$ in sites, $[1\ldots i] \;,
  [i+1\ldots j]$ and $[j+1\ldots n]$ will, respectively, create
  subsequences orderisomorph to forbidden patterns $_f 4231$, $_f
  2431$ and $_f 2413$. Now, from Property \ref{p1h3} $(ii)$, and
  Definition \ref{dseH3} if $j=n$, all permutations resulting from
  insertion of $n+1$ in all but the last active sites will be
  labelized $[B,k]$. With $k$ taking value in $[3\dots t+1]$. Indeed,
  all sites before the insertion are preserved and sites, $j+1$, $n+1$
  are always active (Property \ref{p1h3} (iii)),
\item site $n+1$: from Property \ref{p1h3} (ii), the resulting permutation has label $[A,t+1]$.
\end{itemize}
\item if $\pi$ is labeled $[B,t]$, it has $t$ active sites. Now, if $n+1$ is inserted in:
\begin{itemize}
\item the first $(t-2)^{th}$ sites: again, we apply Property
  \ref{p1h3} $(iii)$ and consequently, these permutations have label
  $[B,k]$. With $k$ taking value in $[3\dots t+1]$,
\item site $j+1$: by Property \ref{p1h3} $(i)$, its label will be
  $[A,1]$,
\item site $n+1$: the resulting permutation has label
  $[A,t]$. Indeed, from Property \ref{p1h3} $(ii)$, the $t-2$ active
  sites before $n$ will be preserved and the last two sites are always
  active.
\end{itemize}
\end{itemize}
\end{proof}

\subsubsection{Generating tree of $S_n({\cal H}_4)$}

\begin{definition}
\label{ssh4}
\label{dseH4}
We consider the following succession system and labeling:

$
\left\{ \begin{array}{lll}
root = [A,2] \\
{[A,t]} & \rightarrow & [B,3],\ldots,[B,t+1],[A,t+1] \\
{[B,t]} & \rightarrow & [B,3],\ldots,[B,t],[A,t],[A,1] \\
\end{array} \right. \\
$
Given $\pi$ in $S_n({\cal H}_4)$ and $i,j$ such that $\pi(i)=n-1$ and
$\pi(j)=n$, $\pi$ has label:
\begin{itemize}
\item $[A,t]$ if $i<j$ and $\pi$ has $t$ active sites,
\item $[B,t]$ if $i>j$ and $\pi$ has $t$ active sites.
\end{itemize}
\end{definition}

\begin{property}
\label{p1H4}
Set $\pi$ in $S_n({\cal H}_4)$:
\begin{itemize}
\item[(i)] if $\pi$ has label $[A,t]$ then $i+1 \leq j \leq i+2$,
\item[(ii)] if $\pi$ has label $[A,t]$ and $j=i+1$. Then all sites
  from the first to the ${j+1}^{th}$ are active, all others being
  inactive,
\item[(iii)] if $\pi$ has label $[A,t]$ and $j=i+2$. Then the only active site is $i+1$,
\item[(iv)] if $\pi$ has label $[B,t]$. All sites from $1$ to $j+2$
  are active. Moreover, all sites from $j+3$ to $n$ are inactive.
\end{itemize} 
\end{property}

\begin{proof} $ $
\begin{itemize}
\item[(i)] set $j>i+2$, the following permutation:
  $(n)\pi(i+1)\pi(i+2)(n+1)$ will be either orderisomorph to $_f 3124$
  or $_f 3214$,
\item[(ii)] consider the following permutation, satisfying the hypothesis:
  $\pi(1)\ldots(n-1)(n)\ldots \pi(n)$. All sites from $j+2$ to $n+1$
  are inactive as the subsequence $\pi(i)\pi(j)\pi(j+1)(n+1)$ is
  orderisomorph to the forbidden pattern $_f 2314$. Now, all sites
  from the first to the ${j+1}^{th}$ are active. As before, we prove
  this by contradiction. If some subsequence $a\, b\, c\, (n+1)$ is
  orderisomorph to forbidden patterns who end with $4$, then, $a\, b\,
  c\, \pi(j)$ would have been orderisomorph to the same patterns. Now,
  if the subsequence $a\, b\, (n+1)\, c$ is orderisomorph to $_f
  2143$. Then, either $a\, b\,\pi(i)\pi(j)$ or $a\, b\, \pi(i)c$ would
  have been too. The remaining case is if a subsequence $(n+1)\, a\,
  b\, c$ is orderisomorph to $_f 4213$, then either, $a\, b\,
  c\,\pi(i)$ or $\pi(i)\,a\, b\, c$ would have been orderisomorph to $_f
  2134$ or $_f 4213$,
\item[(iii)] consider the following permutation, satisfying the
 hypothesis: $\pi(1)\ldots(n-1)\;\pi(i+1)(n)\ldots \pi(n)$. Insertion of $n+1$
  in sites $[1\ldots i]$, $i+2$ and $[j+1\ldots n+1]$ will result with
  subsequences, respectively orderisomorph to $_f 4213$, $_f 2143$ and
  $_f 2134$. Consequently all those sites are inactive. Now, insertion
  of $n+1$ in position $i+1$ cannot result in a subsequence $a\, b\,
  c\, (n+1)$ orderisomorph to $_f 2134$, $_f 2314$, $_f 3124$ nor $_f
  3214$. Indeed, if such a subsequence exists then $a\, b\, c\,
  \pi(j)$ would also have been orderisomorph to the same patterns
  which contradicts the hypothesis. Now, if $a\, b\, (n+1)\, c$ is
  orderisomorph to $_f 2143$ then $a\, b\,\pi(i)\pi(j)$ would have
  been orderisomorph to $_f 2134$. Finally, if $(n+1)\, a\, b\, c$ is
  orderisomorph to $_f 4213$, then $\pi(i)\, a\, b\, c$ would have
  been too, as $c$ cannot correspond to $\pi(j)$ in the subsequence,
\item[(iv)] if $n+1$ is inserted in position $k>j+2$, then the
  following subsequence: $(n)\pi(j+1)\pi(j+2)(n+1)$ will be
  orderisomorph to either $_f 3124$ or $_f 3214$ since all sites
  greater than $j+2$ are inactive. If site $j+1$ is inactive, then
  some subsequence $b\, a\, (n)(n+1)$, with $a<b$, must be
  orderisomorph to $_f 2134$. This is contradiction as
  $b\,a\,(n)(n-1)$ would have been orderisomorph to $_f 2143$. Now, if
  a site $1 \leq k \leq j$ is inactive, then some subsequence
  involving $n$ and $n+1$ must be oderisomorph to either $_f 2143$ or
  $_f 4213$. In both case, we get a contradiction. Finally, if site
  $i+2$ is inactive, then some subsequence $a \, (n)\,\pi(j+1)(n+1)$
  must be orderisomorph to $_f 2314$. Hence a contradiction since
  subsequence $a\, \pi(j+1)\pi(j+2)(n-1)$ would have been
  orderisomorph in $\pi$ to either $_f 2134$, $_f 3124$ or $_f 3214$.
\end{itemize}  
\end{proof}
Now, we prove that given a labelized permutation, its children will have the labels provided by the generating tree.
\begin{proof}
Set $\pi$ in $S_n({\cal H}_4)$:

\begin{itemize}
\item if $\pi$ is labelized $[A,t]$, it has $t$ active sites. Now, if
  $n+1$ is inserted in:
\begin{itemize}
\item the first to the last but one active site $k$: from Definition
  \ref{dseH4} along with Property \ref{p1H4} $(iv)$,the new
  permutation has label $[B,k+2]$,
\item the last site, which is the only active site if $t=1$: from
  Property \ref{p1H4} $(iv)$, the new permutation has label
  $[A,t+1]$.
\end{itemize}
\item if $\pi$ is labelized $[B,t]$, it has $t$ active sites. Now, if
  $n+1$ is inserted in:
\begin{itemize}
\item the first to $(t-2)^{th}$ active site $k$: from  Property
  \ref{p1H4} $(iv)$, the new permutation has label $[B,k+2]$,
\item the last but one active site: from Property \ref{p1H4} $(ii)$,
  the new permutation has label $[A,t+1]$,
\item the last active site: from Property \ref{p1H4} $(iii)$, the
  new permutation has label $[A,1]$.
\end{itemize}
\end{itemize}
\end{proof}

\subsection{Permutations avoiding ${\cal H}_3^{-1\,c}$ and
  ${\cal H}_5$ are characterized by the same succession system} 

\begin{proof}[Proof of Lemma \ref{lss3}] Both succession systems
  associated with the generating trees defined hereafter are
  equivalent up to an isomorphism. As those generating trees produces
  $S_n({\cal H}_3^{-1\,c})$ and $S_n({\cal H}_5)$, both are
  characterized by the same succession system.
\end{proof}

\subsubsection{Generating tree of $S_n({\cal H}_3^{-1\,c})$}

Remember that ${\cal H}_3^{-1\,c}$ is the $inverse$ and $complement$
of ${\cal H}_3$. So we have \\ 
${\cal H}_3^{-1\,c}=\{1324,2134,2143,2314,2413,3214\}$.

\begin{definition}
\label{ssh3ic}
We consider the following succession system and labeling:

$
\left\{ \begin{array}{lll}
root = [A,2] \\
{[P]} & \rightarrow  & [A,2] \\
{[A,t]} & \rightarrow & [B,t+1],[A,3],\ldots,[A,t+1] \\
{[B,t]} & \rightarrow & [B,3],[A,3],[P]^{t-2} \\
\end{array} \right. \\
$
Given $\pi$ in $S_n({\cal H}_3^{-1\,c})$ and integers $i,j$ such that
$\pi(i)=n-1$ and $\pi(j)=n$, $\pi$ is labeled:
\begin{itemize}
\item $[P]$ if $i=1$ and $2<j$,
\item $[A,t]$ if either $2<j<i$ or $j=i+1$,
\item $[B,t]$ if $j=1$.
\end{itemize}
\end{definition}

\begin{property}
\label{pH3ic}
Set $\pi$ in $S_n({\cal H}_3^{-1\,c})$:
\begin{itemize}
\item[(i)] if a site $k$ is inactive, all sites to the right of $k$
  are inactive. Therefore, if $\pi$ has $t$ active sites, all sites
  from the first to $t$ are inactive. All others being inactive,
\item[(ii)] if $j=1$ and $i>j+1$, then all, but the first, sites are
  inactive. As a consequence, if $\pi$ has label $[P]$ only its
  first site is active,
\item[(iii)] if $\pi$ has label $[A,t]$ then site $j+1$ is active and
  site $j+2$ is inactive. Thus, with property $(i)$ above, all sites
  from the first to the $(j+1)^{th}$ are active, all others being
  inactive.  
\end{itemize}
\end{property}

\begin{proof} $ $
\begin{itemize}
\item[(i)] results directly from the structure of the forbidden
  subsequences in ${\cal H}_3^{-1\,c}$,
\item[(ii)] set $\pi$ with label $[P]$: $(n-1)\pi(2)\ldots(n)\ldots
  \pi(n)$. Clearly, all sites from the second to the last will yield a
  subsequence orderisomorph to either: $_f 2413$, $_f 2143$ or $_f
  2134$,
\item[(iii)] site $j+2$ is inactive since either $2<j<i$ and
  then, subsequence $\pi(1)(n)\pi(j+1)(n+1)$ is orderisomorph to $_f
  1324$ or $_f 2314$, or $j=i+1$ and then subsequence
  $(n-1)(n)\pi(j+1)(n+1)$ is orderisomorph to $_f 2314$. Finally, site
  $j+1$ is active, as if not then some subsequence $b\,a\,(n)(n+1)$
  must be orderisomorph to $_f 2134$, from which follows a contradiction as then
  either subsequences, $b\,a\,(n)(n-1)$ or $b\,a\,(n-1)(n)$ would have been
  orderisomorph to $_f 2143$ or $_f 2134$ respectively.
\end{itemize}
\end{proof}

\begin{figure}[h]
\center
\input{Schema_H3ic}
\caption{generating tree of $S_n({\cal H}_3^{-1\,c})$.}
\label{H3ic}
\end{figure}

Now, we prove that given a labelized permutation, its children will have the labels provided by the generating tree.
\begin{proof}
Set $\pi$ in $S_n({\cal H}_3^{-1\,c})$:

\begin{itemize}
\item if $\pi$ is labelized $[P]$. We conclude from Definition
  \ref{ssh3ic} and Properties \ref{pH3ic} $(ii)$ and $(iii)$, that
  the solechild has label $[A,2]$,
\item if $\pi$ is labelized $[A,t]$, it has $t$ active sites. Now, if
  $n+1$ is inserted in:
\begin{itemize}
\item the first site: the resulting permutation has label $[B,t+1]$,
  since all active sites of $\pi$ remain active. Indeed, all sites
  before $j$ necessary remain active. Site $j+1$ too, as no
  subsequence involving $(n+1)(n)(n+2)$ can be orderisomorph to $_f
  1324$ or $_f 3214$,
\item site $k$ in $[2\ldots t]$: from Property \ref{pH3ic} $(iii)$,
  the resulting permutation has label $[A,k+1]$.
\end{itemize}
\item if $\pi$ is labelized $[B,t]$, it has $t$ active sites. Now, if
  $n+1$ is inserted in:
\begin{itemize}
\item the first site: then the resulting permutation has label
  $[B,3]$. Indeed, the permutation will be:
  $(n+1)(n)\pi(2)\ldots\pi(n)$. The first three sites are active and
  the fourth will result with a subsequence, namely,
  $(n+1)(n)\pi(2)(n+2)$, orderisomorph to $_f 3214$,
\item the second site: then the resulting permutation has label $[A,3]$,
\item site $k$ in $[3\ldots t]$: then, from Property \ref{pH3ic}
  $(i)$, the resulting permutation has label $[P]$.
\end{itemize}
\end{itemize}
\end{proof}

\subsubsection{Generating tree of $S_n({\cal H}_5)$}

\begin{definition}
\label{ssh5}
We consider the following succession system and labeling:

$
\left\{ \begin{array}{lll}
root = [A,2] \\
{[P]} & \rightarrow  & [A,2] \\
{[A,t]} & \rightarrow & [A,t+1],[B,t+1],[A,3],\ldots,[A,t] \\
{[B,t]} & \rightarrow & [3,A],[3,B],[P]^{t-2} \\
\end{array} \right. \\
$
Given $\pi$ in $S_n({\cal H}_5)$ and integers $i,j$ such that
$\pi(i)=n-1$ $\pi(j)=n$, $\pi$ is labeled:
\begin{itemize}
\item $[P]$ if $1<i<j$,
\item $[A,t]$ if either, $j \not = 2$ and $1<i<j$ or, $i=1$ and $j=3$,
\item $[B,t]$ if $j=2$.
\end{itemize}
\end{definition}

\begin{property}
\label{p1H5}
Consider $\pi$ in $S_n({\cal H}_5)$, if site $k$ is inactive, then all
sites from $k$ to $n+1$ are inactive. Consequently, if $\pi$ has $t$
active sites, then all sites from the first to the $k^{th}$ are
active.
\end{property}
\begin{proof}
results directly from the structure of the forbidden subsequences in ${\cal H}_5$,
\end{proof}
Now, we proove that given a labelized permutation, its children will have the labels provided by the generating tree.
\begin{proof}
Set $\pi$ in $S_n({\cal H}_5)$:

\begin{itemize}
\item if, $\pi$ has label $[P]$, the first site is active as no forbidden
  patterns begin with $4$. Moreover, the second site is inactive
  as $\pi(1)(n+1)\ldots (n-1) \ldots (n)$ is orderisomorph to $_f 1423$. Next, we consider the only resulting
  permutation: $(n+1) \pi(1) \ldots (n-1) \ldots (n)$. The third site is inactive as it corresponds to the former second
  site. The second site is active as no forbidden patterns begin with $34$. Finally, since 
  $(n+1)$ is in first position, the new label is $[A,2]$.
\item if $\pi$ has label $[A,t]$, it has $t$ active sites. Now, if $n+1$ is inserted in:
\begin{itemize}
\item the first site: the resulting permutation as label
 [A,t+1]. Indeed, if site $t+1$ becomes inactive then some subsequence
 $(n+1)\, a\, b\, (n+2)$ with $(n+2)$ in position $t+1$ and $a<b$ must
 be orderisomorph to $_f 3124$. Hence a contradiction as the
 subsequence formed with $a\,b\,(n+1)$ (with $n+1$ in position $t$) 
and $n$ would have been orderisomorph with either $_f 3124$, $_f 1324$,
 $_f 1234$ or $_f 1243$. Site $t+2$ is inactive as it corresponds to the former site $t+1$ in $\pi$,
\item the second site: the resulting permutation has label
  $[B,t+1]$. As above, if site $t+1$ becomes inactive then some
  subsequence $a\, (n+1)\, b\, (n+2)$ with $(n+2)$ in position $t+1$
  and $a<b$ must be orderisomorph to $_f 1324$ and from this follows a contradiction,
\item site $k$ in $[3\ldots t]$: the resulting permutation has label
  $[A,k]$, as $n+1$ isn't in the first or third place. Moreover site
  $k$ remains active and site $k+1$ is inactive as patterns $_f 1234$
  and $_f 2134$ are forbidden.
\end{itemize}
\item if $\pi$ has label $[B,t]$, it has $t$ active sites. Now, if
  $n+1$ is inserted in:
\begin{itemize}
\item the first site: the resulting permutation,
  $(n+1)\pi(1)(n)\ldots \pi(n)$, has label $[A,3]$, as the third site
  is active and the fourth inactive ($(n+1)\pi(1)(n)(n+2)$ is
  orderisomorph to pattern $_f 3124$),
\item the second site: the resulting permutation,
  $\pi(1)(n+1)(n)\ldots \pi(n)$, has label $[A,3]$ since the third
  site is active and the fourth inactive ($\pi(1)(n+1)(n)(n+2)$ is
  orderisomorph to pattern $_f 1324$),
\item all other active sites: the resulting permutations have label
  $[P]$ as one can find $1<i<j$ such that $\pi(i)=n$ and $\pi(j)=n+1$.
\end{itemize}
\end{itemize}
\end{proof}



\input{bibliographie}
\end{document}

%% file: Commd_bibl.tex
\newcommand{\bibauteur}[1]{{\bf #1}}
\newcommand{\bibtitre}[1]{{#1}}
\newcommand{\bibjournal}[1]{{\it #1}}
\newcommand{\bibvolume}[1]{{\bf #1}}

\newcommand{\bibannee}[1]{{(#1)}}
\newcommand{\bibpages}[2]{{#1}--{#2}}

\newcommand{\DM}{Discrete Mathematics}

\newcommand{\FPSAC}{Formal Power Series and Algebraic Combinatorics}
\newcommand{\emeFPSAC}[1]{Proceedings of $#1$th \FPSAC}

\newcommand{\TCS}{Theorical Computer Science}

\newcommand{\JCT}{Journal of Combinatorial Theory}
\newcommand{\JCTA}{\JCT\ (Series A)}

\newcommand{\BURO}{Bureau Universitaire de Recherche Op\'erationnelle}
\newcommand{\CBURO}{Cahiers du \BURO}

\newcommand{\EurJC}{European Journal of Combinatorics}

\newcommand{\UnivBX}{University Bordeaux~1, France}

\newcommand{\TheseBX}{PHD-thesis, \UnivBX}

\newcommand{\AcadBelge}{M\'emoires publi\'ees par l'Acad\'emie royale de Belgique}
\newcommand{\ElecJC}{The Electronic Journal of Combinatorics}

%% file: Schema_F3.tex

\begin{center}

\begin{picture}(0,85)(65,-55)

\node[linecolor=White,Nh=3.0,Nmr=1.5](n01)(0,-2.5){$1\;[2]$}

\node[linecolor=White,Nh=3.0,Nmr=1.5,Nadjust=w](n02)(40,-25){$\;1\;2\;[3]$}

\node[linecolor=White,Nh=3.0,Nmr=1.5,Nadjust=w](n04)(80,-10){$\;3\;1.2\;[T]$}
\node[linecolor=White,Nh=3.0,Nmr=1.5,Nadjust=w](n05)(80,-25){$\;1\;3.2\;[3]$}
\node[linecolor=White,Nh=3.0,Nmr=1.5,Nadjust=w](n06)(80,-42.5){$\;1\;2\;3\;[4]$}

\node[linecolor=White,Nh=3.0,Nmr=1.5,Nadjust=w](n07)(120,-5){$\;4\;3.1.2\;[T]$}
\node[linecolor=White,Nh=3.0,Nmr=1.5,Nadjust=w](n08)(120,-10){$\;3.4.1.2\;[2]$}
\node[linecolor=White,Nh=3.0,Nmr=1.5,Nadjust=w](n09)(120,-15){$\;3.1.2\;4\;[3]$}
\node[linecolor=White,Nh=3.0,Nmr=1.5,Nadjust=w](n10)(120,-20){$\;4\;1.3.2\;[T]$}
\node[linecolor=White,Nh=3.0,Nmr=1.5,Nadjust=w](n11)(120,-25){$\;1\;4.3.2\;[3]$}
\node[linecolor=White,Nh=3.0,Nmr=1.5,Nadjust=w](n12)(120,-30){$\;1\;3.2\;4\;[4]$}
\node[linecolor=White,Nh=3.0,Nmr=1.5,Nadjust=w](n13)(120,-35){$\;4\;1.2.3\;[T]$}
\node[linecolor=White,Nh=3.0,Nmr=1.5,Nadjust=w](n14)(120,-40){$\;1\;4.2.3\;[3]$}
\node[linecolor=White,Nh=3.0,Nmr=1.5,Nadjust=w](n15)(120,-45){$\;1\;2\;4.3\;[4]$}
\node[linecolor=White,Nh=3.0,Nmr=1.5,Nadjust=w](n16)(120,-50){$\;1\;2\;3\;4\;[5]$}

\node[linecolor=White,Nh=3.0,Nmr=1.5,Nadjust=w](n03)(40,17.5){$\;2\;1\;[T]$}

\node[linecolor=White,Nh=3.0,Nmr=1.5,Nadjust=w](n17)(80,30){$\;3\;2.1\;[T]$}
\node[linecolor=White,Nh=3.0,Nmr=1.5,Nadjust=w](n18)(80,17.5){$\;2.3.1\;[2]$}
\node[linecolor=White,Nh=3.0,Nmr=1.5,Nadjust=w](n19)(80,5){$\;2.1\;3\;[3]$}

\node[linecolor=White,Nh=3.0,Nmr=1.5,Nadjust=w](n20)(120,35){$\;4\;3.2.1\;[T]$}
\node[linecolor=White,Nh=3.0,Nmr=1.5,Nadjust=w](n21)(120,30){$\;3.4.2.1\;[2]$}
\node[linecolor=White,Nh=3.0,Nmr=1.5,Nadjust=w](n22)(120,25){$\;3.2.1\;4\;[3]$}
\node[linecolor=White,Nh=3.0,Nmr=1.5,Nadjust=w](n23)(120,20){$\;4\;2.3.1\;[T]$}
\node[linecolor=White,Nh=3.0,Nmr=1.5,Nadjust=w](n24)(120,15){$\;2.3.1\;4\;[3]$}
\node[linecolor=White,Nh=3.0,Nmr=1.5,Nadjust=w](n25)(120,10){$\;4\;2.1.3\;[T]$}
\node[linecolor=White,Nh=3.0,Nmr=1.5,Nadjust=w](n26)(120,05){$\;2.1\;4.3\;[3]$}
\node[linecolor=White,Nh=3.0,Nmr=1.5,Nadjust=w](n27)(120,00){$\;2.1\;3\;4\;[4]$}

\gasset{AHdist=1.4,AHangle=10.12,AHLength=1.8,AHlength=0.8}
\drawedge[sxo=6,syo=-2,exo=-5.0](n01,n02){}
\drawedge[sxo=6,syo=2,exo=-5.0](n01,n03){}

\drawedge[sxo=5,exo=-5.0](n02,n04){}
\drawedge[sxo=5,exo=-5.0](n02,n05){}
\drawedge[sxo=5,exo=-5.0](n02,n06){}

\drawedge[sxo=5,exo=-5.0](n04,n07){}
\drawedge[sxo=5,exo=-5.0](n04,n08){}
\drawedge[sxo=5,exo=-5.0](n04,n09){}

\drawedge[sxo=5,exo=-5.0](n05,n10){}
\drawedge[sxo=5,exo=-5.0](n05,n11){}
\drawedge[sxo=5,exo=-5.0](n05,n12){}

\drawedge[sxo=5,exo=-5.0](n06,n13){}
\drawedge[sxo=5,exo=-5.0](n06,n14){}
\drawedge[sxo=5,exo=-5.0](n06,n15){}
\drawedge[sxo=5,exo=-5.0](n06,n16){}

\drawedge[sxo=5,exo=-5.0](n03,n17){}
\drawedge[sxo=5,exo=-5.0](n03,n18){}
\drawedge[sxo=5,exo=-5.0](n03,n19){}

\drawedge[sxo=5,exo=-5.0](n17,n20){}
\drawedge[sxo=5,exo=-5.0](n17,n21){}
\drawedge[sxo=5,exo=-5.0](n17,n22){}

\drawedge[sxo=5,exo=-5.0](n18,n23){}
\drawedge[sxo=5,exo=-5.0](n18,n24){}

\drawedge[sxo=5,exo=-5.0](n19,n25){}
\drawedge[sxo=5,exo=-5.0](n19,n26){}
\drawedge[sxo=5,exo=-5.0](n19,n27){}

\end{picture}
\end{center}

%% file: Schema_preuve.tex
\begin{center}
\begin{picture}(81,48)(0,-48)

\node[linecolor=Black,linewidth=0.40,NLangle=0.0,Nmr=0.0](n0)(12.0,-12.0){${\cal H}_1$}

\node[linecolor=Black,linewidth=0.40,NLangle=0.0,Nmr=0.0](n3)(44.0,-36.0){${\cal H}_4$}

\node[linecolor=Black,linewidth=0.40,NLangle=0.0,Nmr=0.0](n4)(44.0,-28.0){${\cal H}_3$}

\node[linecolor=Black,linewidth=0.40,NLangle=0.0,Nmr=0.0](n5)(44.0,-20.0){${\cal H}_2$}

\node[linecolor=Black,linewidth=0.40,NLangle=0.0,Nmr=0.0](n13)(44.0,-12.0){${\cal H}_1^{\star}$}

\node[linecolor=Black,linewidth=0.40,NLangle=0.0,Nmr=0.0](n14)(76.0,-28.0){${\cal H}_3^{\mbox{-}1c}$}

\node[linecolor=Black,linewidth=0.40,NLangle=0.0,Nmr=0.0](n18)(76.0,-36.0){${\cal H}_5$}

\drawedge[linewidth=0.20,ELdist=2.6,AHlength=0.45](n4,n14){$-1\,c$}
\drawedge[linewidth=0.20,ELdist=2.6,AHlength=0.45](n14,n4){}

\drawedge[linewidth=0.20,ELdist=2.3,AHlength=0.4](n0,n13){$\star$}
\drawedge[linewidth=0.20,ELdist=2.3,AHlength=0.4](n13,n0){}

\node[dash={2.0 2.0 2.0 3.0}{0.0},ExtNL=y,NLangle=90.0,NLdist=2.0,Nw=16.0,Nh=16.0,Nmr=4.0](n28)(12.0,-12.0){Lemma \ref{lss1}}

\node[dash={2.0 2.0 2.0 3.0}{0.0},ExtNL=y,NLangle=90.0,NLdist=2.0,Nw=16.0,Nh=40.0,Nmr=4.0](n33)(44.0,-24.0){Lemma \ref{lss2}}

\node[dash={2.0 2.0 2.0 3.0}{0.0},ExtNL=y,NLangle=90.0,NLdist=2.0,Nw=16.0,Nh=24.0,Nmr=4.0](n35)(76.0,-32.0){Lemma \ref{lss3}}

\end{picture}
\end{center}

%% file: Schema_H1star.tex

\begin{center}

\begin{picture}(0,85)(65,-55)

\node[linecolor=White,Nh=3.0,Nmr=1.5](n01)(0,-2.5){$1\;[A,2]$}

\node[linecolor=White,Nh=3.0,Nmr=1.5,Nadjust=w](n02)(40,-25){$\;1\;2\;[A,3]$}

\node[linecolor=White,Nh=3.0,Nmr=1.5,Nadjust=w](n04)(80,-10){$\;3\;1.2\;[B,3]$}
\node[linecolor=White,Nh=3.0,Nmr=1.5,Nadjust=w](n05)(80,-25){$\;1\;3\;2\;[B,4]$}
\node[linecolor=White,Nh=3.0,Nmr=1.5,Nadjust=w](n06)(80,-40){$\;1\;2\;3\;[A,4]$}

\node[linecolor=White,Nh=3.0,Nmr=1.5,Nadjust=w](n07)(120,0){$\;4\;3.1.2\;[B,3]$}
\node[linecolor=White,Nh=3.0,Nmr=1.5,Nadjust=w](n08)(120,-5){$.3.4.1.2\;[A,1]$}
\node[linecolor=White,Nh=3.0,Nmr=1.5,Nadjust=w](n09)(120,-10){$\;3.1.2\;4\;[A,3]$}
\node[linecolor=White,Nh=3.0,Nmr=1.5,Nadjust=w](n10)(120,-15){$\;4\;1.3.2\;[B,3]$}
\node[linecolor=White,Nh=3.0,Nmr=1.5,Nadjust=w](n11)(120,-20){$\;1\;4\;3.2\;[B,4]$}
\node[linecolor=White,Nh=3.0,Nmr=1.5,Nadjust=w](n12)(120,-25){$.1.3.4.2\;[A,1]$}
\node[linecolor=White,Nh=3.0,Nmr=1.5,Nadjust=w](n13)(120,-30){$\;1\;3.2\;4\;[A,4]$}
\node[linecolor=White,Nh=3.0,Nmr=1.5,Nadjust=w](n14)(120,-35){$\;4\;1.2.3\;[B,3]$}
\node[linecolor=White,Nh=3.0,Nmr=1.5,Nadjust=w](n15)(120,-40){$\;1\;4\;2.3\;[B,4]$}
\node[linecolor=White,Nh=3.0,Nmr=1.5,Nadjust=w](n16)(120,-45){$\;1\;2\;4\;3\;[B,5]$}
\node[linecolor=White,Nh=3.0,Nmr=1.5,Nadjust=w](n17)(120,-50){$\;1\;2\;3\;4\;[A,5]$}

\node[linecolor=White,Nh=3.0,Nmr=1.5,Nadjust=w](n03)(40,20){$\;2\;1\;[B,3]$}

\node[linecolor=White,Nh=3.0,Nmr=1.5,Nadjust=w](n18)(80,30){$\;3\;2.1\;[B,3]$}
\node[linecolor=White,Nh=3.0,Nmr=1.5,Nadjust=w](n19)(80,20){$.2.3.1\;[A,1]$}
\node[linecolor=White,Nh=3.0,Nmr=1.5,Nadjust=w](n20)(80,10){$\;2.1\;3\;[A,3]$}

\node[linecolor=White,Nh=3.0,Nmr=1.5,Nadjust=w](n21)(120,35){$\;4\;3.2.1\;[B,3]$}
\node[linecolor=White,Nh=3.0,Nmr=1.5,Nadjust=w](n22)(120,30){$.3.4.2.1\;[A,1]$}
\node[linecolor=White,Nh=3.0,Nmr=1.5,Nadjust=w](n23)(120,25){$\;3.2.1\;4\;[A,3]$}
\node[linecolor=White,Nh=3.0,Nmr=1.5,Nadjust=w](n24)(120,20){$.2.3.1\;4\;[A,2]$}
\node[linecolor=White,Nh=3.0,Nmr=1.5,Nadjust=w](n25)(120,15){$\;4\;2.1.3\;[B,3]$}
\node[linecolor=White,Nh=3.0,Nmr=1.5,Nadjust=w](n26)(120,10){$\;2.1\;4\;3\;[B,4]$}
\node[linecolor=White,Nh=3.0,Nmr=1.5,Nadjust=w](n27)(120,05){$\;2.1\;3\;4\;[A,4]$}

\gasset{AHdist=1.4,AHangle=10.12,AHLength=1.8,AHlength=0.8}
\drawedge[sxo=6,syo=-2,exo=-5.0](n01,n02){}
\drawedge[sxo=6,syo=2,exo=-5.0](n01,n03){}

\drawedge[sxo=5,exo=-5.0](n02,n04){}
\drawedge[sxo=5,exo=-5.0](n02,n05){}
\drawedge[sxo=5,exo=-5.0](n02,n06){}

\drawedge[sxo=5,exo=-5.0](n04,n07){}
\drawedge[sxo=5,exo=-5.0](n04,n08){}
\drawedge[sxo=5,exo=-5.0](n04,n09){}

\drawedge[sxo=5,exo=-5.0](n05,n10){}
\drawedge[sxo=5,exo=-5.0](n05,n11){}
\drawedge[sxo=5,exo=-5.0](n05,n12){}
\drawedge[sxo=5,exo=-5.0](n05,n13){}

\drawedge[sxo=5,exo=-5.0](n06,n14){}
\drawedge[sxo=5,exo=-5.0](n06,n15){}
\drawedge[sxo=5,exo=-5.0](n06,n16){}
\drawedge[sxo=5,exo=-5.0](n06,n17){}

\drawedge[sxo=5,exo=-5.0](n03,n18){}
\drawedge[sxo=5,exo=-5.0](n03,n19){}
\drawedge[sxo=5,exo=-5.0](n03,n20){}

\drawedge[sxo=5,exo=-5.0](n18,n21){}
\drawedge[sxo=5,exo=-5.0](n18,n22){}
\drawedge[sxo=5,exo=-5.0](n18,n23){}

\drawedge[sxo=5,exo=-5.0](n19,n24){}

\drawedge[sxo=5,exo=-5.0](n20,n25){}
\drawedge[sxo=5,exo=-5.0](n20,n26){}
\drawedge[sxo=5,exo=-5.0](n20,n27){}

\end{picture}
\end{center}

%% file: Schema_H3ic.tex

\begin{center}

\begin{picture}(0,85)(65,-55)

\node[linecolor=White,Nh=3.0,Nmr=1.5](n01)(0,-2.5){$1\;[A,2]$}

\node[linecolor=White,Nh=3.0,Nmr=1.5,Nadjust=w](n02)(40,-25){$\;1\;2\;[A,3]$}

\node[linecolor=White,Nh=3.0,Nmr=1.5,Nadjust=w](n04)(80,-10){$\;3\;1\;2\;[B,4]$}
\node[linecolor=White,Nh=3.0,Nmr=1.5,Nadjust=w](n05)(80,-25){$\;1\;3\;2.[A,3]$}
\node[linecolor=White,Nh=3.0,Nmr=1.5,Nadjust=w](n06)(80,-40){$\;1\;2\;3\;[A,4]$}

\node[linecolor=White,Nh=3.0,Nmr=1.5,Nadjust=w](n07)(120,0){$\;4\;3\;1.2.[B,3]$}
\node[linecolor=White,Nh=3.0,Nmr=1.5,Nadjust=w](n08)(120,-5){$\;3\;4\;1.2.[A,3]$}
\node[linecolor=White,Nh=3.0,Nmr=1.5,Nadjust=w](n09)(120,-10){$\;3.1.4.2.[P]\;\;\;\,$}
\node[linecolor=White,Nh=3.0,Nmr=1.5,Nadjust=w](n10)(120,-15){$\;3.1.2.4.[P]\;\;\;\,$}
\node[linecolor=White,Nh=3.0,Nmr=1.5,Nadjust=w](n11)(120,-20){$\;4\;1\;3\;2.[B,4]$}
\node[linecolor=White,Nh=3.0,Nmr=1.5,Nadjust=w](n12)(120,-25){$\;1\;4\;3.2.[A,3]$}
\node[linecolor=White,Nh=3.0,Nmr=1.5,Nadjust=w](n13)(120,-30){$\;1\;3\;4\;2.[A,4]$}
\node[linecolor=White,Nh=3.0,Nmr=1.5,Nadjust=w](n14)(120,-35){$\;4\;1\;2\;3\;[B,5]$}
\node[linecolor=White,Nh=3.0,Nmr=1.5,Nadjust=w](n15)(120,-40){$\;1\;4\;2.3.[A,3]$}
\node[linecolor=White,Nh=3.0,Nmr=1.5,Nadjust=w](n16)(120,-45){$\;1\;2\;4\;3.[A,4]$}
\node[linecolor=White,Nh=3.0,Nmr=1.5,Nadjust=w](n17)(120,-50){$\;1\;2\;3\;4\;[B,5]$}

\node[linecolor=White,Nh=3.0,Nmr=1.5,Nadjust=w](n03)(40,20){$\;2\;1\;[B,3]$}

\node[linecolor=White,Nh=3.0,Nmr=1.5,Nadjust=w](n18)(80,30){$\;3\;2\;1.[B,3]$}
\node[linecolor=White,Nh=3.0,Nmr=1.5,Nadjust=w](n19)(80,20){$\;2\;3\;1.[A,3]$}
\node[linecolor=White,Nh=3.0,Nmr=1.5,Nadjust=w](n20)(80,10){$\;2.1.3.[P]\;\;\;\,$}

\node[linecolor=White,Nh=3.0,Nmr=1.5,Nadjust=w](n21)(120,35){$\;4\;3\;2.1.[B,3]$}
\node[linecolor=White,Nh=3.0,Nmr=1.5,Nadjust=w](n22)(120,30){$\;3\;4\;2.1.[A,3]$}
\node[linecolor=White,Nh=3.0,Nmr=1.5,Nadjust=w](n23)(120,25){$\;3.2.4.1.[P]\;\;\;\,$}
\node[linecolor=White,Nh=3.0,Nmr=1.5,Nadjust=w](n24)(120,20){$\;4\;2\;3\;1.[B,4]$}
\node[linecolor=White,Nh=3.0,Nmr=1.5,Nadjust=w](n25)(120,15){$\;2\;4\;3.1.[A,3]$}
\node[linecolor=White,Nh=3.0,Nmr=1.5,Nadjust=w](n26)(120,10){$\;2\;3\;4\;1.[A,4]$}
\node[linecolor=White,Nh=3.0,Nmr=1.5,Nadjust=w](n27)(120,05){$\;4\;2.1.3.[B,2]$}

\gasset{AHdist=1.4,AHangle=10.12,AHLength=1.8,AHlength=0.8}
\drawedge[sxo=6,syo=-2,exo=-5.0](n01,n02){}
\drawedge[sxo=6,syo=2,exo=-5.0](n01,n03){}

\drawedge[sxo=5,exo=-5.0](n02,n04){}
\drawedge[sxo=5,exo=-5.0](n02,n05){}
\drawedge[sxo=5,exo=-5.0](n02,n06){}

\drawedge[sxo=5,exo=-5.0](n04,n07){}
\drawedge[sxo=5,exo=-5.0](n04,n08){}
\drawedge[sxo=5,exo=-5.0](n04,n09){}
\drawedge[sxo=5,exo=-5.0](n04,n10){}

\drawedge[sxo=5,exo=-5.0](n05,n11){}
\drawedge[sxo=5,exo=-5.0](n05,n12){}
\drawedge[sxo=5,exo=-5.0](n05,n13){}

\drawedge[sxo=5,exo=-5.0](n06,n14){}
\drawedge[sxo=5,exo=-5.0](n06,n15){}
\drawedge[sxo=5,exo=-5.0](n06,n16){}
\drawedge[sxo=5,exo=-5.0](n06,n17){}

\drawedge[sxo=5,exo=-5.0](n03,n18){}
\drawedge[sxo=5,exo=-5.0](n03,n19){}
\drawedge[sxo=5,exo=-5.0](n03,n20){}

\drawedge[sxo=5,exo=-5.0](n18,n21){}
\drawedge[sxo=5,exo=-5.0](n18,n22){}
\drawedge[sxo=5,exo=-5.0](n18,n23){}

\drawedge[sxo=5,exo=-5.0](n19,n24){}
\drawedge[sxo=5,exo=-5.0](n19,n25){}
\drawedge[sxo=5,exo=-5.0](n19,n26){}

\drawedge[sxo=5,exo=-5.0](n20,n27){}

\end{picture}
\end{center}

%% file: FineRapport.bbl
\begin{thebibliography}{}

\bibitem
	[BB-MDFGG-B]
	{Banderier}
	\bibauteur{C.~Banderier} 
		\bibauteur{M.~Bousquet--M\'elou}, 
		\bibauteur{A.~Denise}, 
		\bibauteur{P.~Flajolet}, 
		\bibauteur{D.~Gardy} and 
		\bibauteur{D.~Gouyou--Beauchamps},
	\bibtitre{On generating functions of generating trees},
	\bibjournal{\emeFPSAC{11}},
	Barcelona, Spain 
	\bibannee{1999}
	\bibpages{40}{52}.

\bibitem
	[BdLP]
	{Barcucci2}
	\bibauteur{E.~Barcucci}, \bibauteur{A.~Del~Lungo} and 
		\bibauteur{E.~Pergola},
	\bibtitre{Random generation of trees and other 
		combinatorial objects},
	\bibjournal{\TCS}
	\bibvolume{218}
	\bibannee{1999}
	\bibpages{219}{232}.

\bibitem
	[B-M]
	{Bousquet}
	\bibauteur{M.~Bousquet-Mélou},
	\bibtitre{Four classes of pattern-avoiding permutations under
	one roof: generating trees with two labels},
	\bibjournal{\ElecJC}
	\bibvolume{9(2)}
	\bibannee{2003}.


\bibitem
	[BdLPP]
	{Barcucci1}
	\bibauteur{E.~Barcucci} 
		\bibauteur{A.~Del~Lungo}, 
		\bibauteur{E.~Pergola} and 
		\bibauteur{R.~Pinzani},
	\bibtitre{ECO: a methodology for the Enumeration of 
		Combinatorial Objects},
	\bibjournal{Journal of Difference Equations and 
		Applications}
	\bibvolume{5}
	\bibannee{1999}
	\bibpages{435}{490}.

\bibitem
	[CGHK]
	{Chung}
	\bibauteur{F.R.K.~Chung}, \bibauteur{R.L.~Graham}, 
		\bibauteur{V.E.~Hoggat} and 
		\bibauteur{M. Kleiman}, 
	\bibtitre{The number of Baxter permutations}, 
	\bibjournal{\JCTA}
	\bibvolume{24}
	\bibannee{1978}
	\bibpages{382}{394}.

\bibitem
	[C]
	{Callan}
	\bibauteur{D.~Callan},
	\bibtitre{Some identities for the Catalan and Fine numbers},
	\bibvolume{arXiv CO/0502532}
	\bibannee{2005}.

\bibitem
	[DS]
	{Deutsch}
	\bibauteur{E.~Deutsch} and \bibauteur{L.~Shapiro}, 
	\bibtitre{A survey of Fine numbers}, 
	\bibjournal{Discrete Mathematics}	   	
	\bibvolume{241 issue 1-3}
	\bibannee{2001}
	\bibpages{241}{265}.

\bibitem
	[EV]
	{Elder}
	\bibauteur{M.~Elder} and \bibauteur{V.~Vatter},
	\bibtitre{Problems and conjectures presented at The third
	international conference on permutation patterns},
	University of Florida, 
	\bibannee{March 7-11, 2005}.

\bibitem
	[E]
	{Errera}
	\bibauteur{A.~Errera},
	\bibtitre{Un probl\`eme d'\'enum\'eration},
	\bibjournal{\AcadBelge},
	Bruxelles, 
	tome \bibvolume{11}
	\bibannee{1931}.

\bibitem
	[F]
	{Fine}
	\bibauteur{T.~Fine},
	\bibtitre{Extrapolation when very little is known about the source},
	\bibjournal{Information and Control}
	\bibvolume{16}
	\bibannee{1970}
	\bibpages{331}{359}.

\bibitem
	[Gi]
	{GireThese}
	\bibauteur{S.~Gire},
	\bibtitre{Arbres, permutations \`a motifs exclus et cartes 
		planaires~: quelques probl\`emes algorithmiques 
		et combinatoires},
	\bibjournal{\TheseBX}
	\bibannee{1993}.

\bibitem
	[Gu]
	{GuibertThese}
	\bibauteur{O.~Guibert},
	\bibtitre{Combinatoire des permutations \`a motifs exclus 
		en liaison avec mots, cartes planaires et tableaux 
		de Young},
	\bibjournal{\TheseBX}
	\bibannee{1995}.

\bibitem
	[Krew]
	{KrewEventail}
	\bibauteur{G.~Kreweras},
	\bibtitre{Sur les \'eventails de segments},
	\bibjournal{\CBURO}
	\bibvolume{15}
	\bibannee{1970}
	\bibpages{1}{41}.

\bibitem
	[P]
	{PergolaThese}
	\bibauteur{E.~Pergola},
	\bibtitre{ECO: a Methodology for Enumerating Combinatorial Objects},
	\bibjournal{PHD-thesis, University of Firenze, Italy}
	\bibannee{1998}.

\bibitem
	[Ro]
	{Rogers}
	\bibauteur{D.G.~Rogers},
	\bibtitre{Similarity Relations on Finite Ordered Sets},
	\bibjournal{\JCTA}
	\bibvolume{23}
	\bibannee{1977}
	\bibpages{88}{98}.

\bibitem
	[Sh]
	{Shapiro}
	\bibauteur{L.W.~Shapiro},
	\bibtitre{A Catalan triangle},
	\bibjournal{\DM}
	\bibvolume{14}
	\bibannee{1976}
	\bibpages{83}{90}.

\bibitem
	[Si]
	{Simion}
	\bibauteur{R.~Simion} and \bibauteur{F.W.~Schmidt},
	\bibtitre{Restricted permutations},
	\bibjournal{\EurJC}
	\bibvolume{6}
	\bibannee{1985}
	\bibpages{383}{406}.

\bibitem
	[SP]
	{SloanePlouffe}
	\bibauteur{N.J.A.~Sloane} and 
		\bibauteur{S.~Plouffe},
	\bibtitre{The encyclopedia of integer sequences},
	Academic Press
	\bibannee{1995}.

\bibitem
	[Str]
	{Strehl}
	\bibauteur{V.~Strehl},
	\bibtitre{A note on Similarity Relations},
	\bibjournal{\DM}
	\bibvolume{19}
	\bibannee{1977}
	\bibpages{99}{101}.

\bibitem
	[W1]
	{WestPHD}
	\bibauteur{J.~West},
	\bibtitre{Permutations with forbidden subsequences
		and stack-sortable permutations},
	PHD-thesis, Massachusetts Institute of Technology, 
		Cambridge
	\bibannee{1990}.

\bibitem
	[W2]
	{WestCatalan}
	\bibauteur{J.~West},
	\bibtitre{Generating trees and 
		the Catalan and Schr\"{o}der numbers},
	\bibjournal{\DM}
	\bibvolume{146}
	\bibannee{1995}
	\bibpages{247}{262}.

\bibitem
	[W3]
	{WestAG}
	\bibauteur{J.~West},
	\bibtitre{Generating trees and forbidden subsequences},
	\bibjournal{\DM}
	\bibvolume{157}
	\bibannee{1996}
	\bibpages{363}{374}.

\end{thebibliography}
